\documentclass{amsart}
\usepackage[english]{babel}
\usepackage[utf8]{inputenc}
\usepackage[T1]{fontenc}

\usepackage{a4wide}             
\usepackage{emptypage}          
\usepackage{amsmath}            
\usepackage{amssymb}            
\usepackage{mathtools}          
\allowdisplaybreaks             
\usepackage{amsthm}             
\usepackage{graphicx}           
\usepackage{tikz}               
\usetikzlibrary{arrows.meta}               
\usetikzlibrary{decorations.pathreplacing} 
\usetikzlibrary{decorations.markings}
\usepackage{dsfont}             
\usepackage{float}              
\usepackage{esint}              
\usepackage{comment}            
\usepackage{etoolbox}           

\usepackage{hyperref}
\hypersetup{
    colorlinks=true,
    linkcolor=black,
    filecolor=magenta,      
    urlcolor=cyan,
    citecolor=black
}
\usepackage[noabbrev]{cleveref} 

\numberwithin{equation}{section}

\theoremstyle{plain}
\newtheorem{theorem}{Theorem}[section]
\newtheorem{lemma}[theorem]{Lemma}
\newtheorem{proposition}[theorem]{Proposition}
\newtheorem{corollary}[theorem]{Corollary}

\theoremstyle{definition}
\newtheorem{definition}[theorem]{Definition}
\newtheorem{remark}[theorem]{Remark}
\newtheorem{hypothesis}{Condition}

\AtBeginEnvironment{lemma}{\crefalias{theorem}{lemma}}
\AtBeginEnvironment{proposition}{\crefalias{theorem}{proposition}}
\AtBeginEnvironment{corollary}{\crefalias{theorem}{corollary}}
\AtBeginEnvironment{definition}{\crefalias{theorem}{definition}}
\AtBeginEnvironment{remark}{\crefalias{theorem}{remark}}

\crefname{hypothesis}{Condition}{Conditions}
\Crefname{hypothesis}{Condition}{Conditions}

\AfterEndEnvironment{definition}{\noindent\ignorespaces}    
\AfterEndEnvironment{proposition}{\noindent\ignorespaces}
\AfterEndEnvironment{corollary}{\noindent\ignorespaces}
\AfterEndEnvironment{remark}{\noindent\ignorespaces}
\AfterEndEnvironment{lemma}{\noindent\ignorespaces}
\AfterEndEnvironment{theorem}{\noindent\ignorespaces}
\AfterEndEnvironment{hypothesis}{\noindent\ignorespaces}
\AfterEndEnvironment{proof}{\noindent\ignorespaces}

\makeatletter
\newcommand{\sectionnotoc}[1]{\begingroup\let\@tocwrite\@gobbletwo\section*{#1}\endgroup}
\makeatother

\newcommand{\virgolette}[1]{``#1''}

\newcommand{\R}{\mathbb{R}}
\newcommand{\N}{\mathbb{N}}
\newcommand{\Z}{\mathbb{Z}}
\newcommand{\ud}{\mathrm{d}}

\newcommand{\Hh}{\mathbb{H}}

\newcommand{\suchthat}{\;\ifnum\currentgrouptype=16 \middle\fi|\;}

\newcommand{\eps}{\varepsilon}

\DeclareMathOperator{\curl}{curl}
\DeclareMathOperator{\Curl}{Curl}

\DeclareMathOperator{\dist}{dist}

\title[Energetic scaling in Grain Boundaries via Dirichlet-to-Neumann map]{Non-local Interfacial Energies in Grain Boundaries via Dirichlet-to-Neumann Map and Energetic scaling}
\author{Valerio Corica, Francesco Santilli}
\begin{document}

\begin{abstract}
This paper deals with the variational analysis of edge dislocations arising in small-angle grain boundaries (Read--Shockley) and semi-coherent interfaces (van der Merwe). Using the Nabarro--Peierls approach to model edge dislocations, we reduce the two-dimensional isotropic Navier--Lamé system governing the linearized elastic energy to a one-dimensional non-local energy on the slip interface through a Dirichlet-to-Neumann map. We prove that the resulting boundary energy coincides, up to a multiplicative constant, with the fractional \(H^{\frac12}\)-Gagliardo seminorm of the displacement jump along the interface. Within this framework, we establish the existence of a unique global minimizer by combining our dimensional-reduction result with the characterization of De Luca, Goldman, and Ponsiglione. Indeed, we show that the optimal defect configuration is uniquely attained, up to translations, by an equispaced periodic array. Finally, we provide a physical rescaling that recovers the classical logarithmic energy regimes of both the Read--Shockley and van der Merwe models.
\par\vspace{0.3cm}
\par\vspace{0.1cm}
\end{abstract}

\maketitle

\tableofcontents

\section*{Introduction}
The study of crystal defects, and in particular dislocations, is a central theme in materials science and solid mechanics, as these line defects represent the primary carriers of plastic deformation in crystalline solids \cite{Hirth83}. While classical continuum elasticity provides an accurate description of the stress and strain fields away from the defects, it is no longer valid close to the dislocation lines; indeed, it predicts diverging elastic energies at the so-called defect cores. To bridge the gap between the discrete atomistic structure of crystals and continuum elastic theories, the seminal models of Nabarro and Peierls \cite{Nabarro47,Peierls40} introduced a multiscale description in which the crystal is treated as a linear elastic continuum away from a prescribed slip plane, while a non-linear potential accounts for the atomic mismatch across the interface.
In this paper, we aim to describe the emergence and the energetic characterization of dislocations in two distinct physical configurations: \textit{grain boundaries} and \textit{semi-coherent interfaces}. Both frameworks can be viewed, in a simplified setting, as the matching of two square crystal lattices along a one-dimensional interface.

For \textit{small-angle grain boundaries}, as described by the classical Read--Shockley theory \cite{Read50}, the two phases consist of identical lattices that are slightly rotated with respect to each other, for instance by angles $\theta$ and $-\theta$. In the case of \textit{semi-coherent interfaces}, whose classical description dates back to van der Merwe \cite{van50}, the lattices have the same orientation but different spacing, resulting in a mismatch. In both cases, the geometric incompatibility between the two bulk lattices prevents a coherent matching across the interface. The accumulation of macroscopic elastic strain is then relieved through the formation of arrays of dislocations localized close to the interface.

From a variational point of view, this mechanism can be described through the competition between bulk elastic energy and the energetic cost associated to the formation of defects. A first rigorous derivation of the Read--Shockley logarithmic scaling without any \textit{a priori} prescription on the number, position, or nature of the defects was obtained in the seminal work of Lauteri and Luckhaus \cite{Lauteri16}. Their result is composed of a delicate lower bound estimate and an explicit upper bound of the scaling, which has also been obtained by Scardia and Tolotti \cite{Scardia26} by means of a different construction. Building on their deep analysis, Fortuna, Garroni, and Spadaro \cite{Garroni25} subsequently upgraded their work to a \(\Gamma\)-convergence result for \textit{grain boundaries}. This sharp interface approach has also inspired recent phase-field models designed to incorporate crystallographic symmetries \cite{Conti26}. In a parallel direction, the first author, together with E. Spadaro, is currently studying a Lauteri--Luckhaus-type model for semi-coherent interfaces, with the aim of deriving the corresponding van der Merwe scaling. The resulting picture is consistent with the analysis developed in the present paper and suggests closely related qualitative mechanisms for \textit{grain boundaries} and \textit{semi-coherent interfaces}. Recently, the Read--Shockley scaling was also recovered in \cite{Grabner25} by minimizing the discrete elastic energy associated to the Volterra potential.

A complementary mathematical viewpoint arises from the dimensional reduction of bulk elastic problems to non-local energies supported on the interface. In the context of straight parallel dislocations in a three-dimensional crystal, translation invariance along the dislocation lines reduces the problem to a two-dimensional one, set on a cross-section orthogonal to them. We will consider a cylindrical domain of the type \(\Hh_T \coloneqq \Hh^+_T \cup \Hh^-_T\), where \( \Hh^\pm_T \coloneqq [0,T] \times \R^\pm\) with \(T \in \N\) fixed, and in each half-plane there is a crystal associated to a displacement field \( U \coloneqq U^{+} \mathds{1}_{\Hh^+_T} + U^{-} \mathds{1}_{\Hh^-_T} \). In our model, the elastic equilibrium is governed by the Navier--Lamé equations. The bulk elastic energy of a displacement field can be expressed solely in terms of the jump of the trace on the interface $w \coloneqq u^+-u^- \in H^{\frac{1}{2}}_T(\R)$, where \(u^\pm (x)= U^\pm(x,0)\) in the sense of traces, through the associated Dirichlet-to-Neumann map. The resulting one-dimensional energy is proportional to the fractional $H^{\frac{1}{2}}$-Gagliardo seminorm of the interfacial slip $w$. As suggested in \cite{Garroni05} in the periodic setting and proved in a non-periodic framework \cite{cinesi20}, this reduction provides a rigorous connection between Nabarro--Peierls-type models and non-local interfacial ones.

On the other hand, the origin of the interfacial defects can be understood through the compatibility constraints imposed by the two lattices. As shown in \cite{Fanzon20}, a small geometric mismatch gives rise to a topological constraint on the admissible interfacial displacement \(U\). Such a constraint forces the accumulation of a prescribed phase mismatch, which can be accommodated by a fixed number of deltas modeling dislocations; in our model, the \(T\)-periodicity assumption forces the function \(u\) to have \(T\) deltas in the periodicity interval \([0,T)\).

Dimension reduction results, such as the so-called Dirichlet-to-Neumann map with respect to the extension problem of the Laplace equation, are well-known tools \cite{Caffarelli07,Palatucci12}; nevertheless, we show the latter in the periodic setting for \(s=\frac{1}{2}\) providing a proof using Fourier series techniques, which appears to be missing in the current literature; indeed, in \cite{Caffarelli07} the result is briefly recalled using functional analysis arguments in a non-periodic setting. Given this fact and inspired by the Fourier transform arguments in \cite{cinesi20}, we adopt once again Fourier series techniques also in the elastic 2D problem, where the elastic energy associated to the Navier--Lamé equations reduces
to a \((\frac{1}{2},2)\)-Gagliardo seminorm. Nonetheless, the map does not require many assumptions on the energy density on the interface: provided a suitable symmetry of the displacement \(U\) on both half-planes \(\Hh^\pm_T\), the result can be extended to a wider class of models. 

Driven also by the aforementioned equivalence and the seminal work of Müller \cite{Muller93}, the variational analysis of the non-local fractional energies and the existence of periodic minimizers has motivated a vast body of work such as \cite{Goldman26,Pini26,DeLuca24,Santilli26,RenWei03}. In particular, while a macroscopic periodicity or a bounded domain is typically prescribed \textit{a priori}, the competition between interaction terms favors the emergence of an intrinsic micro-periodicity. In our context, the state is typically described by a profile \(u\), which can be thought of as a displacement slip. Regarding phase-field models, one can require \(u' \in \{a,b\}\) almost everywhere, with \(a,b \neq 0\). In \cite{DeLuca24}, imposing \(a=1 \text{ and } b=-\Lambda\) with \(\Lambda \gg 1\), the authors characterize the optimal periodicity using variational techniques, while in \cite{Giuliani11}, imposing \(a=1=-b\), the result is achieved via \textit{reflection positivity}, thanks to the symmetry of the slopes. In \cite{Goldman26}, imposing \(u'=1-\sum_{i=1}^T \delta_{x_i}\), the authors study the \(H^s\)-Gagliardo seminorm in the sub-critical and critical regimes \(s \in (0,\frac{1}{2}]\) and manage to prove the uniqueness of the periodic minimizer by studying the first and second variation of the functional.
Relying on this last characterization, we focus our study on the latter constraint and exploit the leading order of the energy of the minimizer \(u=x-\lfloor x \rfloor\) (up to translations). Then, using a proper scaling, we are able to recover the Read--Shockley formula for \textit{small-angle grain boundaries} and the van der Merwe theory for \textit{semi-coherent interfaces}, as already obtained for the latter in \cite{DeLuca24}. 

In summary, the main contribution of this work is to provide a unified variational bridge between two-dimensional linear elasticity and one-dimensional fractional interfacial energies in a periodic setting. By developing the Dirichlet-to-Neumann reduction on the strip $\Hh_T$ via discrete Fourier series, our framework simultaneously accommodates normal interfacial jumps (for Read--Shockley tilt boundaries) and tangential jumps (for van der Merwe interfaces). Combined with the uniqueness of periodic equispaced ground states, this yields a derivation of the classical Read--Shockley logarithmic scaling $\theta|\log\theta|$.

\par\vspace{5pt}\noindent
\textbf{Structure of the paper.} The paper is organized as follows.
In \Cref{sec: section1}, we set up the geometric framework, collect the periodic function spaces, and establish the Dirichlet-to-Neumann map for the Laplacian on the periodic strip.
In \Cref{sec: section2}, we present the 2D linear elastic model and establish the Dirichlet-to-Neumann dimensional reduction.
In \Cref{sec: section3}, we investigate the non-local minimization problem on the interface and the uniqueness of equispaced configurations.
Finally, in \Cref{sec: section4}, we perform the energetic rescaling and recover the celebrated Read--Shockley and van der Merwe results.

\section{Notation and Setting} \label{sec: section1} \label{sec: notation}

\begin{figure}[htpb]
\centering
\begin{tikzpicture}[
	scale=1.2,
	defect/.style={circle, fill=red, inner sep=1.8pt},
	interface/.style={very thick, black}
	]
	\fill[blue!5] (0,0) rectangle (6,3);
	\fill[red!5]  (0,-3) rectangle (6,0);
	
	\fill[blue!2] (-1.5,0) rectangle (0,3);
	\fill[blue!2] (6,0) rectangle (7.5,3);
	\fill[red!2] (-1.5,-3) rectangle (0,0);
	\fill[red!2] (6,-3) rectangle (7.5,0);
	
	\draw[interface] (-1.5, 0) -- (7.5, 0) node[right] {$\Hh^0_T$};
	
	\draw[dashed, thick, darkgray] (0, -3.2) -- (0, 3.2) node[above, black] {$x=0$};
	\draw[dashed, thick, darkgray] (6, -3.2) -- (6, 3.2) node[above, black] {$x=T$};
	
	\draw[decorate, decoration={brace, amplitude=8pt, raise=4pt}, thick]
	(0, 3.5) -- (6, 3.5) node [black, midway, yshift=18pt] {Periodic strip $\Hh_T$};
	
	\node[defect, label=below:{$x_1$}] at (1,0) {};
	\node[defect, label=below:{$x_2$}] at (3.5,0) {};
	\node[defect, label=below:{$x_T$}] at (5.2,0) {};
	
	\node[defect, opacity=0.3] at (-1.5,0) {};
	\node[defect, opacity=0.3] at (-0.8,0) {}; 
	\node[defect, opacity=0.3] at (7,0) {};   
	
	\node[font=\Large] at (3, 1.5) {$\Hh^+_T$};
	\node[font=\Large] at (3, -1.5) {$\Hh^-_T$};
	
	\draw[-{Latex[length=2mm]}, thick] (-1, -2.5) -- (-1, -1.5) node[above] {$y$};
	\draw[-{Latex[length=2mm]}, thick] (-1, -2.5) -- (0, -2.5) node[right] {$x$};
	
\end{tikzpicture}
\caption{Schematic representation of the periodic setting. The defects $x_i$ lie on the straight interface $\Hh^0_T$, which separates the upper and lower periodic half-planes $\Hh^+_T$ and $\Hh^-_T$.}
\end{figure}
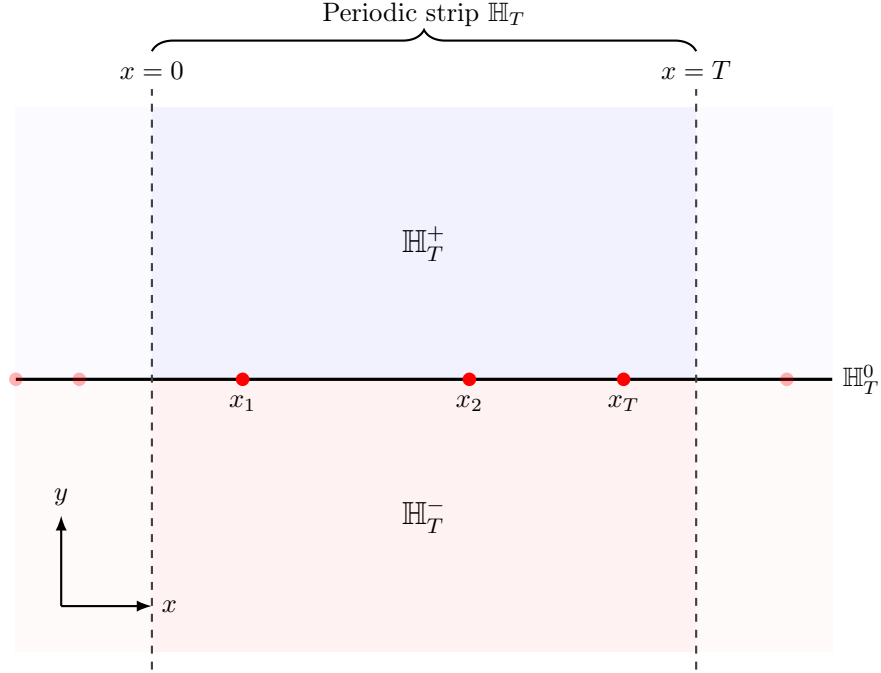

For the reader's convenience, the main notation used in this work is collected in this section. Given \(T>0\), we denote by \(\fint u = \frac{1}{T}\int_0^T u(x) \ud x\) the average of a function \(u\) over the periodic cell. Moreover, we define the periodic semi-strips \(\Hh_T^\pm \coloneqq [0,T] \times \R^\pm\) and their interface \(\Hh^0_T \coloneqq [0,T] \times \{0\} \cong [0,T]\). 

We define the space of \(T\)-periodic \(L^2\) functions on \(\R\) as
\begin{equation*}
L^2_{T}(\R)\coloneqq \left\{u:\R\to\R \text{ measurable} \suchthat u(x+T)=u(x) \text{ for a.e. } x\in\R, \text{ and } u \in L^2(0,T) \right\},
\end{equation*}
endowed with the standard norm \(\|u\|_{L^2(0,T)}\). Furthermore, we define the periodic Sobolev space \(H^1_T(\Hh^\pm_T)\) as
\begin{equation*}
	H^1_T(\Hh^\pm_T) \coloneqq \left\{ U \in H^1_{\mathrm{loc}}(\R \times \R^\pm) \suchthat 
	\begin{aligned}
		& U(x+T, y) = U(x, y) \text{ for a.e. } (x,y) \in \R \times \R^\pm, \\
		& U \in L^2(\Hh^\pm_T), \ \nabla U \in L^2(\Hh^\pm_T, \R^2)
	\end{aligned}
	\right\}.
\end{equation*}
Consequently, we define the homogeneous periodic Sobolev space \(\dot{H}^1_T(\Hh^\pm_T)\) as
\begin{equation*}
	\dot{H}^1_T(\Hh^\pm_T) \coloneqq \left\{ U \in H^1_{\mathrm{loc}}(\R \times \R^\pm) \suchthat 
	\begin{aligned}
		& U(x+T, y) = U(x, y) \text{ for a.e. } (x,y) \in \R \times \R^\pm, \\
		& \nabla U \in L^2(\Hh^\pm_T, \R^2)
	\end{aligned}
	\right\}.
\end{equation*}
Similarly, we introduce the \(T\)-periodic fractional Sobolev space
\begin{equation*}
H^{\frac{1}{2}}_{T}(\R)\coloneqq \left\{ u\in L^2_{T}(\R) \suchthat [u]_{\frac{1}{2},T} <+\infty \right\},
\end{equation*}
where \([u]_{\frac{1}{2},T}\) denotes the periodic Gagliardo seminorm:
\begin{equation*}
[u]_{\frac{1}{2},T} \coloneqq \left(\int_0^T \int_{\R} \frac{|u(x)-u(y)|^2}{|x-y|^{2}} \, \ud y \ud x\right)^{\frac{1}{2}},
\end{equation*}
endowed with the standard norm \(\|u\|_{H^{\frac{1}{2}}_T} \coloneqq \|u\|_{L^2(0,T)} + [u]_{\frac{1}{2},T}\). Furthermore, we define $H^{-\frac{1}{2}}_T(\R)$ as the topological dual space of $H^{\frac{1}{2}}_T(\R)$. Its duality pairing $\langle \cdot, \cdot \rangle_{H^{-\frac{1}{2}}_T \times H^{\frac{1}{2}}_T}$ can be identified with the continuous extension of the $L^2$ scalar product, i.e., for every $v \in H^{-\frac{1}{2}}_T(\R)$ and $u \in H^{\frac{1}{2}}_T(\R)$:
\begin{equation*}
    |\langle v, u \rangle| \le \|v\|_{H^{-\frac{1}{2}}_T} \|u\|_{H^{\frac{1}{2}}_T} < +\infty,
\end{equation*}
where $\langle v, u \rangle = \int_0^T v(x)u(x) \, \ud x$ whenever $v \in L^2_T(\R)$. The space $H^{-\frac{1}{2}}_T(\R)$ is endowed with the standard dual norm.

Let us consider \(g \in H^{\frac{1}{2}}_{T}(\Hh^0_T)\). We denote the space of \(T\)-periodic \(H^1\)-functions with boundary datum \(g\) on \(\Hh^0_T\) as 
\begin{equation*}
H_{g,T}^1(\Hh^\pm_T) \coloneqq \left\{ u \in H^1_T(\Hh^\pm_T) \suchthat u|_{\Hh^0_T}=g \right\},
\end{equation*}
and similarly we denote by \(\dot{H}^1_{g,T}(\Hh^\pm_T)\) the corresponding space for homogeneous Sobolev functions. Moreover, we define the space of periodic harmonic functions as 
\begin{equation*}
A(\Hh^\pm_T) \coloneqq \{ V \in H^1_T(\Hh^\pm_T) \mid \Delta V = 0\}.
\end{equation*}
As is well known, the \(\curl\) of a distribution \(T \in \mathcal{D}'(\Omega, \R^2)\) for a set \(\Omega \subset \R^2\) is defined as
\begin{equation*}
	\langle \curl T, \varphi \rangle_{\mathcal{D}' \times \mathcal{D}} \coloneqq - \langle T, \nabla^{\perp} \varphi \rangle_{\mathcal{D}' \times \mathcal{D}}\,,
\end{equation*} 
for every \(\varphi \in C^{\infty}_c(\Omega,\R)\).
We will also consider this operator for matrix fields, and in that case it is understood row-wise and denoted by \(\Curl\).
\begin{remark}
	It is worth noting that, if \(\Omega \subset \R^2\) is a simply connected domain, the row-wise curl-free condition for a matrix-valued distribution \(T \in \mathcal{D}'(\Omega, \R^{2 \times 2})\) is equivalent to \(T\) being a distributional gradient. That is, \(\Curl T = 0\) in \(\mathcal{D}' (\Omega, \R^2)\) if and only if there exists a vector-valued distribution \(U \in \mathcal{D}'(\Omega, \R^2)\) such that \(T = D U\) in the sense of distributions.
\end{remark}
\begin{remark}
	Let \(\Omega \subset \R^2\) and let \(F \in L^2(\Omega, \R^2)\). Then \(\curl F \in H^{-1}(\Omega)\). This follows from the definition of distributional \(\curl\) and the Cauchy--Schwarz inequality; indeed, for any \(\varphi \in C^{\infty}_c(\Omega)\), we have
	\begin{equation*}
		|\langle \curl F , \varphi \rangle_{\mathcal{D}'\times \mathcal{D}}| \leq \int_{\Omega} |F| |\nabla^\perp \varphi| \, \ud x \leq \|F\|_{L^2(\Omega)} \|\varphi\|_{H_0^1(\Omega)}.
	\end{equation*}
	By the density of \(C^{\infty}_c(\Omega)\) in \(H^1_0(\Omega)\), this functional extends continuously to the whole \(H^1_0(\Omega)\). In particular, this allows us to regard $\curl F$ as an $H^{-1}$ object.
\end{remark}

\subsection{Preliminary results}
In this subsection, some fundamental results useful throughout the whole paper are presented.
\begin{lemma} \label{lemma: normal_trace_harmonic}
	Let \(V \in A(\Hh^+_T)\). 
	For every \(\psi \in H^{\frac{1}{2}}_T(\Hh^0_T)\) and for every \(\varphi \in H^1_{\psi,T}(\Hh^+_T)\), the linear operator \(V \mapsto \partial_y V(\cdot, 0)\), defined by the identity
	\begin{equation} \label{eq: trace_harmonic}
		\langle \partial_y V(\cdot, 0), \psi \rangle_{H^{-\frac{1}{2}}_T \times H^{\frac{1}{2}}_T} \coloneqq - \int_{\Hh^+_T} \nabla V \cdot \nabla \varphi \ud x \ud y,
	\end{equation}
	is continuous from \( A(\Hh^+_T)\) to \(H^{-\frac{1}{2}}_T(\Hh^0_T)\).
\end{lemma}

\begin{proof}
	The operator is well-posed: let \(\psi \in H^{\frac{1}{2}}_T(\Hh^0_T)\) and consider two arbitrary extensions \(\varphi_1, \varphi_2 \in H^1_{\psi,T}(\Hh^+_T)\). We define their difference \(\eta \coloneqq \varphi_1 - \varphi_2 \in H^1_{0,T}(\Hh^+_T)\).
	Evaluating the integral on the right-hand side of \cref{eq: trace_harmonic} on \(\eta\), we get:
	\[
	- \int_{\Hh^+_T} \nabla V \cdot \nabla \eta \, \ud x \ud y=0,
	\]
	where we used the harmonicity of \(V\). This proves that the right-hand side of \cref{eq: trace_harmonic} does not depend on the extension.
	
		To prove continuity, fix \(\psi \in H^{\frac{1}{2}}_T(\Hh^0_T)\) such that \(\|\psi\|_{H^{\frac{1}{2}}_T(\Hh^0_T)} = 1\). By a standard construction involving a cut-off function in the $y$-variable, one can define a continuous extension operator \(E \colon H^{\frac{1}{2}}_T(\Hh^0_T) \to H^1_T(\Hh^+_T)\). Taking \(\varphi \coloneqq E\psi\), we have \(\varphi \in H^1_{\psi,T}(\Hh^+_T)\) and there exists a constant \(K>0\) such that \(\|\varphi\|_{H^1_T(\Hh^+_T)} \leq K \|\psi\|_{H^{\frac{1}{2}}_T(\Hh^0_T)}\).
	Therefore, using the Cauchy--Schwarz inequality, we obtain:
	\begin{align*}
		|\langle \partial_y V(\cdot, 0), \psi \rangle|
		&\leq \|\nabla V\|_{L^2(\Hh^+_T)} \|\nabla \varphi\|_{L^2(\Hh^+_T)} \\
		&\leq \|V\|_{H^1_T(\Hh^+_T)} \|\varphi\|_{H^1_T(\Hh^+_T)}\\
		&\leq K \|V\|_{H^1_T(\Hh^+_T)}.
	\end{align*}
	Taking the supremum over all \(\psi\), this concludes the proof. 
\end{proof}

Given \(u \in H^{\frac{1}{2}}_{T}(\R)\), we define its \textit{harmonic extension} to the half-plane \(\Hh^\pm_T\) as the unique function \(U \in \dot{H}^1_T(\Hh^\pm_T)\) which solves
\begin{equation} \label{eq: laplaciano_con_dirichlet}
\begin{cases}
	\Delta U (x,y) = 0 \qquad &\text{in } \Hh^\pm_T\,, \\
	U(x,0)=u(x) \qquad &\text{on } \Hh^0_T\,.
\end{cases}
\end{equation}

The second result deals with the Dirichlet-to-Neumann map with respect to the Laplacian operator, already cited in \cite{Caffarelli07,Palatucci12}, which ensures that \(U\) belongs to the homogeneous Sobolev space \(\dot H^1_T(\Hh^\pm_T)\).
\begin{theorem}[Laplacian Dirichlet-to-Neumann map] \label{thm: DtN sulla striscia}
Let \(u \in H^{\frac{1}{2}}_{T}(\R)\) and let \(U\) be its harmonic extension on \(\Hh^\pm_T\). Then \(U \in \dot{H}^1_T(\Hh^\pm_T)\) and the following identity holds true:
\begin{equation} \label{eq: dtn}
	\int_{\Hh^\pm_T} |\nabla U|^2 \ud x \ud y = \frac{1}{2\pi} [u]^2_{\frac{1}{2},T}.
\end{equation}
Moreover, \(U\) minimizes the Dirichlet energy among all valid extensions with zero boundary condition, i.e.,
\begin{equation} \label{eq: pb minimo dirichlet}
	\min_{\substack{V \in \dot{H}^1_{0,T}(\Hh^\pm_T)}}  \int_{\Hh^\pm_T} |\nabla (U+V)|^2 \ud x \ud y = \int_{\Hh^\pm_T} |\nabla U|^2 \ud x \ud y.
\end{equation}
Furthermore, if \(\fint u(x) \ud x = 0\), then \(U \in H^1_T(\Hh^\pm_T)\).
\end{theorem}
\begin{proof}
By symmetry, the proof is given for \(\Hh^+_T\), using an argument which involves Fourier series techniques. To simplify the notation, let us define the angular frequencies \(\omega_n \coloneqq \frac{2n\pi}{T}\). Since \(u \in L^2_T(\R)\), we can express it through its complex Fourier series as:
\begin{equation*}
	u(x) = \sum_{n=-\infty}^{+\infty} \hat{u}_n e^{i\omega_n x},
\end{equation*}
where \(\hat{u}_0 = \fint_0^T u(x) \ud x\) is the average of \(u\) and, more generally, \(\hat{u}_n\) are its Fourier coefficients. To find the harmonic extension \(U(x,y)\) in the upper half-plane \(\Hh^+_T\), it is enough to solve \cref{eq: laplaciano_con_dirichlet} by separation of variables, requiring \(T\)-periodicity in \(x\) and boundedness of the energy as \(y \to +\infty\), obtaining the following unique solution:
\begin{equation*} 
	U(x, y) = \sum_{n=-\infty}^{+\infty} \hat{u}_n e^{-|\omega_n| y} e^{i\omega_n x}.
\end{equation*}
Now, we evaluate the Dirichlet energy of \(U\), using integration by parts on the truncated domain \(A_L\coloneqq [0,T] \times [0, L]\) and, then, taking the limit as \(L \to +\infty\):
\begin{align}
	\nonumber 
	\int_{\Hh^+_T} |\nabla U|^2 \, \ud x\,\ud y 
	&= \lim_{L\to+\infty} \int_{\partial A_L} U \frac{\partial U}{\partial \mathbf{n}}\, \ud \sigma \nonumber\\ 
	&= \lim_{L\to+\infty} \left(
	\int_0^L U(T,y) \frac{\partial U}{\partial x}(T,y)\, \ud y 
	- \int_0^L U(0,y) \frac{\partial U}{\partial x}(0,y)\, \ud y \right. \nonumber\\
	&\qquad\qquad 
	-\left. \int_0^T U(x,0) \frac{\partial U}{\partial y}(x,0)\, \ud x
	+ \int_0^T U(x,L) \frac{\partial U}{\partial y}(x,L)\, \ud x 
	\right) \nonumber\\
	&= - \int_0^T u(x) \frac{\partial U}{\partial y}(x,0)\, \ud x  \label{eq: identita energia}.
\end{align}
In the first equality we have used the harmonicity of \(U\), whereas in the last equality the boundary terms on the lateral sides cancel out due to periodicity and the term at \(y=L\) vanishes due to the exponential decay as \(L \to +\infty\) (see \Cref{fig: domain_AL}). Every integral is meant in the sense \(H^{-\frac{1}{2}}\times H^{\frac{1}{2}}\).
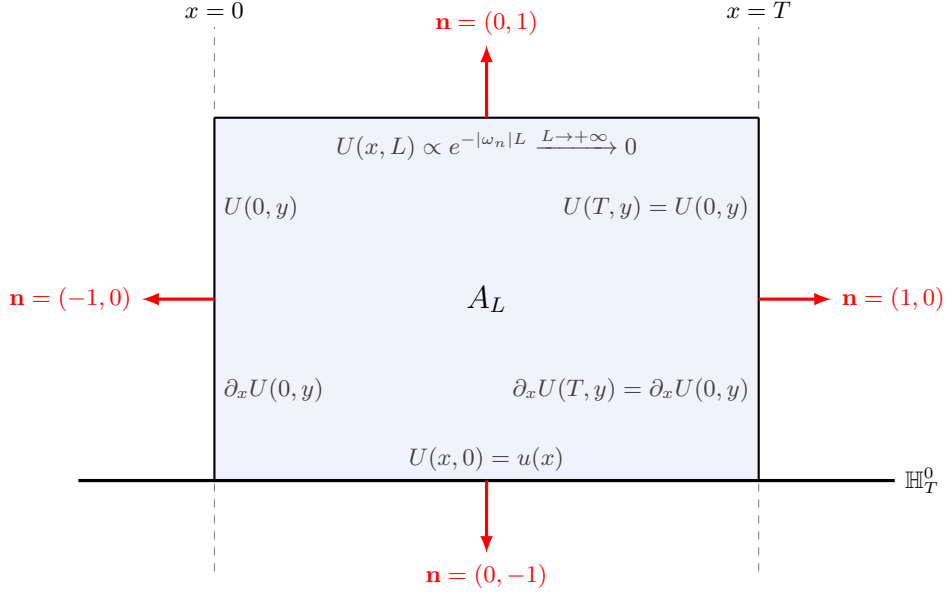
\begin{figure}[ht]
	\centering
	\begin{tikzpicture}[
		scale=1.2,
		normal/.style={->, >=latex, very thick, red},
		boundary/.style={thick, black},
		interface/.style={very thick, black},
		font=\small
		]
		\fill[blue!5] (0,0) rectangle (6,4);
		
		\draw[interface] (-1.5,0) -- (7.5,0) node[right] {$\Hh^0_T$};
		\draw[dashed, gray] (0,-1) -- (0,5) node[above, black] {$x=0$};
		\draw[dashed, gray] (6,-1) -- (6,5) node[above, black] {$x=T$};
		
		\draw[boundary] (6,0) -- (6,4); 
		\draw[boundary] (6,4) -- (0,4); 
		\draw[boundary] (0,4) -- (0,0); 
		
		\node[font=\Large] at (3,2) {$A_L$};
		
		
		\draw[normal] (3,0) -- (3,-0.8) node[below] {$\mathbf{n} = (0,-1)$};
		\node[above, darkgray] at (3, 0) {$U(x,0) = u(x)$};
		
		\draw[normal] (3,4) -- (3,4.8) node[above] {$\mathbf{n} = (0,1)$};
		\node[below, darkgray] at (3, 4) {$U(x,L) \propto e^{-|\omega_n|L} \xrightarrow{L \to +\infty} 0$};
		
		\draw[normal] (0,2) -- (-0.8,2) node[left] {$\mathbf{n} = (-1,0)$};
		\node[right, darkgray, align=left] at (0, 3) {$U(0,y)$};
		\node[right, darkgray, align=left] at (0, 1) {$\partial_x U(0,y)$};
		
		\draw[normal] (6,2) -- (6.8,2) node[right] {$\mathbf{n} = (1,0)$};
		\node[left, darkgray, align=right] at (6, 3) {$U(T,y) = U(0,y)$};
		\node[left, darkgray, align=right] at (6, 1) {$\partial_x U(T,y) = \partial_x U(0,y)$};
		
	\end{tikzpicture}
	\caption{Visual representation of the integration domain $A_L = [0,T] \times [0,L]$ for the Dirichlet energy with the outward normal vectors.}
	\label{fig: domain_AL}
\end{figure}
We now proceed to evaluate the normal derivative of \(U\) at \(y=0\), interchanging the derivative operator and the sum in view of \Cref{lemma: normal_trace_harmonic}:
\begin{equation} \label{eq: derivata normale sul bordo}
	\frac{\partial U}{\partial y}(x, 0) = \sum_{n=-\infty}^{+\infty} (-|\omega_n|) \hat{u}_n e^{i\omega_n x}.
\end{equation}
Substituting \cref{eq: derivata normale sul bordo} into \cref{eq: identita energia} and applying Parseval identity leads to
\begin{align*}
	\int_{\Hh^+_T} |\nabla U|^2 \ud x \ud y &= - \int_0^T u(x) \frac{\partial U}{\partial y} (x, 0) \ud x \\
	&= - T \sum_{n=-\infty}^{+\infty} \hat{u}_n \overline{\left(-|\omega_n| \hat{u}_n\right)} = T \sum_{n=-\infty}^{+\infty} |\omega_n| |\hat{u}_n|^2.
\end{align*}
Next, we want to prove the following identity:
\begin{equation} \label{eq: fourier_seminorm}
	T \sum_{n=-\infty}^{+\infty} |\omega_n| |\hat{u}_n|^2 = \frac{1}{2\pi} [u]^2_{\frac{1}{2},T}.
\end{equation}
In order to do so, we proceed by computing the Fourier series of \(u(x+h)\), which is
\begin{equation*}
	u(x+h) = \sum_{n=-\infty}^{+\infty} \hat{u}_n e^{i\omega_n (x+h)} = \sum_{n=-\infty}^{+\infty} \left( \hat{u}_n e^{i\omega_n h} \right) e^{i\omega_n x}.
\end{equation*}
Since \(\hat{u}_n(e^{i\omega_n h}-1)\) are the Fourier coefficients of the function \( x \mapsto u(x+h)-u(x)\), Parseval identity yields:
\begin{align*}
	\int_0^T |u(x+h)-u(x)|^2 \ud x &= T \sum_{n=-\infty}^{+\infty} |\hat{u}_n|^2 \left| e^{i\omega_n h} - 1 \right|^2 \\
	&= T \sum_{n=-\infty}^{+\infty} |\hat{u}_n|^2 \left[ (\cos(\omega_n h) - 1)^2 + \sin^2(\omega_n h) \right] \\
	&= T \sum_{n=-\infty}^{+\infty} |\hat{u}_n|^2 4\sin^2\left(\frac{\omega_n h}{2}\right).
\end{align*}
Dividing by \(h^2\) and integrating over \(\R\) with respect to \(h\), we obtain:
\begin{equation*} 
	[u]^2_{\frac{1}{2},T} = \int_\R T \sum_{n=-\infty}^{+\infty} |\hat{u}_n|^2 \frac{4 \sin^2 \left( \frac{\omega_n h}{2} \right)}{h^2} \ud h.
\end{equation*}
Exchanging the summation and the integral via Tonelli Theorem and using the fundamental integral
\cite[integral 3.821.9]{Grad14}
\[
\int_{-\infty}^{+\infty} \frac{\sin^2(\alpha h)}{h^2} \ud h = |\alpha|\pi,
\]
which holds for every \(\alpha\in \R\), we obtain:
\begin{equation*}
	[u]^2_{\frac{1}{2},T} = T \sum_{n=-\infty}^{+\infty} |\hat{u}_n|^2 \left( 4 \frac{|\omega_n|}{2} \pi \right) = 2\pi T \sum_{n=-\infty}^{+\infty} |\omega_n| |\hat{u}_n|^2,
\end{equation*}
which implies \cref{eq: fourier_seminorm} and completes the proof of \cref{eq: dtn}. 

Regarding \cref{eq: pb minimo dirichlet}, let \(V \in \dot{H}^1_{0,T}(\Hh^+_T)\). Then, expanding the Dirichlet energy we get:
\begin{align*}
	\int_{\Hh^+_T} |\nabla (U+V)|^2 \ud x \ud y =& \int_{\Hh^+_T} |\nabla U|^2 \ud x \ud y + \int_{\Hh^+_T} |\nabla V|^2 \ud x \ud y + 2 \int_{\Hh^+_T} \nabla U \cdot \nabla V \ud x \ud y\\
	=& \int_{\Hh^+_T} |\nabla U|^2 \ud x \ud y + \int_{\Hh^+_T} |\nabla V|^2 \ud x \ud y \\
	&+ 2\left( - \int_{\Hh^+_T} (\Delta U) V \ud x \ud y + \int_{\partial \Hh^+_T} V \frac{\partial U}{\partial \mathbf{n}} \ud \sigma \right)\\
	=& 	\int_{\Hh^+_T} |\nabla U|^2 \ud x \ud y + \int_{\Hh^+_T} |\nabla V|^2 \ud x \ud y \geq \int_{\Hh^+_T} |\nabla U|^2 \ud x \ud y,
\end{align*}
where in the last equality we have used that \(U\) is harmonic and that the integral on the boundary is zero: indeed, on \(\{y=0\}\), \(V\) vanishes and, on the lateral boundaries \(\{x=0\}\) and \(\{x=T\}\), the function \(U\) is periodic with opposite normal derivative.

To prove the last assertion, assume \(\fint_0^T u(x) \ud x = 0\), which directly implies \(\hat{u}_0 = 0\). Now, we are left to prove that \(\|U\|_{L^2(\Hh^+_T)}\) is finite. Parseval identity on \([0,T]\) yields
\begin{align*}
	\int_0^{+\infty} \int_0^T |U(x,y)|^2 \ud x \ud y &= T \sum_{n \in \mathbb{Z} \setminus \{0\}} |\hat{u}_n|^2 \int_0^{+\infty} e^{-2|\omega_n| y} \ud y \\
	&= \frac{T}{2} \sum_{n \in \mathbb{Z} \setminus \{0\}} |\hat{u}_n|^2 \frac{1}{|\omega_n|} \\
	&= \frac{T^2}{4\pi} \sum_{n \in \mathbb{Z} \setminus \{0\}} \frac{|\hat{u}_n|^2}{|n|} \leq \frac{T^2}{4\pi} \sum_{n \in \mathbb{Z} \setminus \{0\}} |\hat{u}_n|^2 < +\infty,
\end{align*}
where in the last inequality we used the fact that \(u \in H^{\frac{1}{2}}_T(\R) \subset L^2_T(\R)\). This proves that \(U \in H^1_T(\Hh^+_T)\).
\end{proof}
\begin{remark}\label{rmk: h12infourier}
For a generic $T$-periodic function $u \in H^{\frac{1}{2}}_T(\R)$, the following identity holds:
\begin{equation*}  
	[u]_{\frac{1}{2}, T}^2 = 2\pi T \sum_{n \in \mathbb{Z}} |\omega_n| |\hat{u}_n|^2 = 4\pi^2 \sum_{n \in \mathbb{Z}} |n| |\hat{u}_n|^2,
\end{equation*}
where \(\omega_n=\frac{2\pi n}{T}\). In particular, if $u$ has zero average, then $\hat{u}_0 = 0$.
This last result justifies the Fourier-based representation (as in \cite{Ambrosio20}); indeed, the following holds
 \begin{equation*}
 	H^{\frac{1}{2}}_T(\R) = \left\{ u \in L^2_T(\R) \suchthat \sum_{n \in \mathbb{Z}} |n| |\hat{u}_n|^2 < +\infty \right\}.
 \end{equation*}
\end{remark}

\begin{remark} \label{rmk: media nulla}
As shown in the proof, the assumption of zero average of \(u\) guarantees that its harmonic extension \(U\) on the upper and lower half-planes is \(L^2\)-integrable. Moreover, the energy we will deal with is well-posed for every \(U \in \dot{H}^1_T(\Hh^+_T)\).
Nevertheless, we prefer to set the problem in the Hilbert space \(H^1_T(\Hh^+_T)\), so the zero-average condition of the boundary datum will be a fundamental assumption throughout the whole paper.
\end{remark}

\section{The Model} \label{sec: section2} \label{sec: model}
While dislocations in crystals are intrinsically three-dimensional line defects, it is a consolidated approach in the physical and mathematical literature to study them through two-dimensional simplified models by exploiting the symmetries of the defect configuration \cite{Hirth83}. In this paper, we restrict our analysis to straight \textit{edge dislocations}. 
\subsection{Edge dislocations and plane strain reduction}
An edge dislocation can be modeled as an infinite straight line defect in a 3D crystal lattice. By aligning the defect line with the $z$-axis, the system gains translational invariance along this direction. Consequently, the study of the deformations reduces to a 2D \textit{plane strain} problem on a cross-section $\Omega \subset \R^2$ orthogonal to the defect line. In this 2D framework, the classical approach assumes dislocations as a vector-valued empirical measure $\mu$ supported on a set of points $\{x_i\} \subset \Omega$. Namely,
\[\mu = \sum_i \mathbf{b} \delta_{x_i},\]
where $\mathbf{b} \in \R^2$ is the \textit{Burgers vector} characterizing the topological charge of the defects. The presence of such defects induces an elastic strain field $\beta \colon \Omega \to \R^{2 \times 2}$ that must satisfy the incompatibility equation $\Curl \beta = \mu$ in the sense of distributions. A well-known issue in linear elasticity is that the elastic energy associated to a discrete dislocation $\mu$, which is proportional to the $L^2$-norm of the strain $\|\beta\|_{L^2}^2$, diverges logarithmically at the defect cores \cite{Garroni08,Scardia14}. To circumvent this infinite-energy issue and to manage the study of the problem, the standard strategy is to introduce a microscopic core radius $\eps > 0$ (comparable to the lattice spacing) and to replace the singular measure with a mollified continuous defect density $f_\eps \in L^2$, which in our case will be \(f_\eps\coloneqq f \ast \rho^\eps\), with \(\rho^\eps\) a rescaled mollifier.
However, in the specific context of semi-coherent interfaces or small-angle grain boundaries, the total strain does not depend only on isolated defects. The geometric mismatch between the two matching lattices imposes a diffused elastic strain across the interface, which can be represented by the Lebesgue measure $\mathcal{H}^1$. Dislocations, which in our model are modeled as negative Dirac deltas and represent the \textit{plastic} region, then nucleate to relax this macroscopic geometric incompatibility. Therefore, the effective boundary datum of the problem is governed by a distributional derivative of the form
\begin{equation} \label{eq: constraint derivata}
\mathrm{D}u = \mathcal{H}^1 - \sum_{i} \delta_{x_i}.
\end{equation}
This hypothesis can be thought of as the prescription of a boundary datum at infinity, depending on the phenomenon being described.
The competition between the \textit{elastic} and \textit{plastic} regions gives rise to the mantissa-type functions that we analyze in our variational setting (which have already been studied in \cite{Goldman26,DeLuca24}).
As already pointed out, the derivation of the interfacial energy is not strictly related to a density of the form in \cref{eq: constraint derivata}, but remains true for a more general class of defects under the hypotheses of the results in \Cref{sec: model,sec: min}.
\subsection{Physical settings: Grain Boundaries vs. Semi-coherent Interfaces}
We now specialize our setting to the case of an interface separating two different crystalline regions. Let $T > 0$ be the macroscopic periodicity parameter. Let us consider two periodic lattices $\Gamma^+$ and $\Gamma^-$ occupying, respectively, the upper and lower half-planes $\Hh^+_T$ and $\Hh^-_T$ and matching at the straight horizontal interface $\Hh^0_T$.
We can associate the two lattices with the non-singular matrices $A^\pm \in \mathrm{GL}(2)$ such that $\Gamma^\pm= A^\pm( \eps \mathbb{Z}^2)$, where $\eps$ is the atomic distance. Let us define the piecewise constant lattice matrix field $A \coloneqq A^+ \mathds{1}_{\Hh^+_T} + A^- \mathds{1}_{\Hh^-_T}$. From the theory of defects, the dislocation density is encoded by the distributional curl of the lattice configuration. Since $A$ exhibits a jump discontinuity across the horizontal interface $\Hh^0_T$, whose unit tangent vector is $\mathbf{e_1}$, the singular part of its curl is a measure concentrated entirely on the interface. Therefore, we get $\Curl A = (A^+ - A^-) \mathbf{e_1} \, \mathcal{H}^1\llcorner \Hh^0_T$.
This compatibility condition corresponds to the classical \textit{Frank--Bilby equation}, well known in the theory of dislocations (see for example \cite[equation 1]{cinesi26}). 
In our setting, the Burgers vector $\mathbf{b}$ satisfies:
\begin{equation}\label{eq: frank_formula}
\mathbf{b} \propto (A^+ - A^-) \mathbf{e_1}.
\end{equation}
Depending on the choice of the structural matrices $A^\pm$, this geometric framework describes two distinct physical scenarios:

\noindent
\textbf{Small-Angle Grain Boundaries (Read--Shockley model):} 
This case \cite{Read50} occurs when two identical lattices are oriented with an opposite small angle $\alpha$ with respect to the interface. Setting $A^{\pm}= R_{\pm \alpha}$, where $R_\beta \in \mathrm{SO}(2)$ is the rotation matrix of angle \(\beta>0\), we can explicitly compute the resulting Burgers vector from \cref{eq: frank_formula}:
\begin{equation*}
\mathbf{b} \propto (R_{+\alpha} - R_{-\alpha}) \mathbf{e_1} = 2\sin\alpha \, \mathbf{e_2},
\end{equation*}
yielding a Burgers vector \textit{orthogonal} to the interface. 
Observe that choosing equal and opposite angles is due to the linearization considered. Indeed, throughout the paper we work with a linearized energy, with the linearization formally performed around the identity. If the average of the two angles were non-zero, one would instead linearize around a different rotation; by frame indifference, however, this would not affect the problem.
\begin{figure}[ht]
	\centering
	\begin{tikzpicture}[scale=1.2, >=latex]
		\definecolor{ottanio}{RGB}{0, 128, 128}
		\definecolor{ardesia}{RGB}{112, 128, 144}
		
		\def\ang{8}
		\pgfmathsetmacro{\dx}{0.5/tan(\ang)}
		
		\clip (-5.8, -3.4) rectangle (5.0, 3.4);
		
		\begin{scope}
			\clip (-4.3, 0) rectangle (4.3, 3.4);
			\begin{scope}[rotate=\ang, shift={(0, -0.5)}]
				\draw[thick, ardesia] (-8,-2) grid (8,5);
				\foreach \x in {-8,...,8} {
					\foreach \y in {-2,...,5} {
						\filldraw[fill=white, draw=black, thick] (\x, \y) circle (2pt);
					}
				}
			\end{scope}
		\end{scope}
		
		\begin{scope}
			\clip (-4.3, 0) rectangle (4.3, -3.4);
			\begin{scope}[rotate=-\ang, shift={(0, 0.5)}]
				\draw[thick, ardesia] (-8,-5) grid (8,2);
				\foreach \x in {-8,...,8} {
					\foreach \y in {-5,...,2} {
						\filldraw[fill=white, draw=black, thick] (\x, \y) circle (2pt);
					}
				}
			\end{scope}
		\end{scope}
		
		\draw[very thick, black] (-4.9, 0) node[left] {$\Hh^0_T$} -- (4.7, 0);
		
		\node[text=ottanio, font=\huge, fill=white, inner sep=1pt] at (-\dx, 0) {$\boldsymbol{\vdash}$};
		\node[text=ottanio, font=\huge, fill=white, inner sep=1pt] at (\dx, 0) {$\boldsymbol{\vdash}$};
		
		\draw[line width=2pt, violet, -{Latex[length=3.5mm, width=2.5mm]}] 
		(0, 0.4) -- (0, 1.4) node[midway, right=1pt, fill=white, inner sep=1pt] {$\mathbf{b}$};
		
	\end{tikzpicture}
	\caption{Atomistic representation of a small-angle symmetric grain boundary (Read--Shockley model).}
\end{figure}
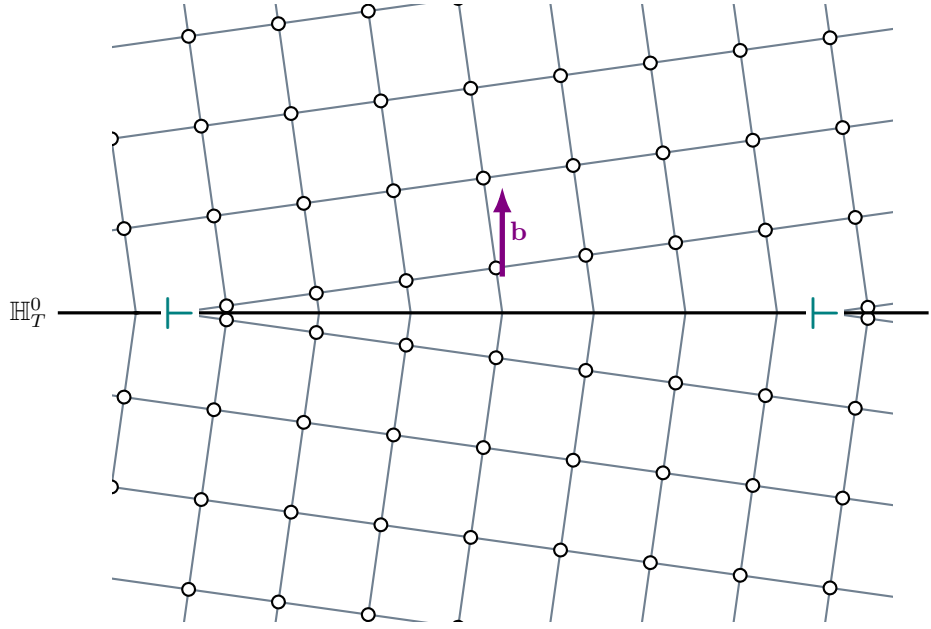

\noindent
\textbf{Semi-coherent Interfaces (van der Merwe model):}
When the two matching lattices share the same spatial orientation but have different lattice parameters, a purely dilatational misfit is generated. Setting $A^+=\operatorname{Id}$ and $A^-= (1+\delta) \operatorname{Id}$ for a small misfit parameter $\delta>0$, we recover the classic semi-coherent interface framework studied in \cite{van50} and recently revisited in a variational setting \cite{DeLuca24}. Applying \cref{eq: frank_formula}, we obtain
\begin{equation*}
\mathbf{b} \propto -\delta \, \mathbf{e_1}.
\end{equation*}
In this setting, the Burgers vector is \textit{parallel} to the interface.
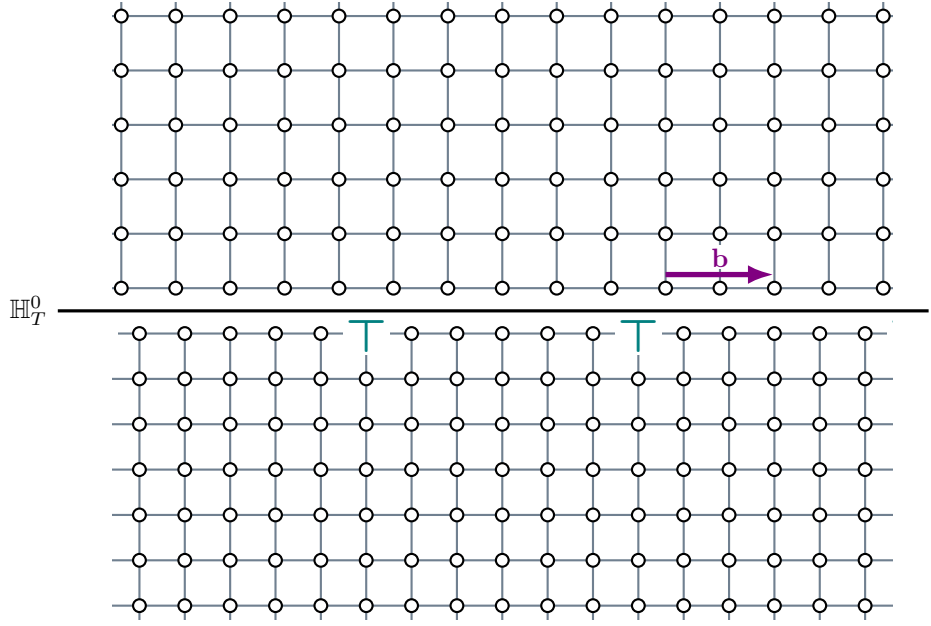
\begin{figure}[ht]
	\centering
	\begin{tikzpicture}[scale=1.2, >=latex]
		\definecolor{ottanio}{RGB}{0, 128, 128}
		\definecolor{ardesia}{RGB}{112, 128, 144}
		
		\clip (-5.8, -3.4) rectangle (5.0, 3.4);
		
		\def\dplus{0.6}   
		\def\dminus{0.5}  
		\def\bound{0.25}  
		
		\begin{scope}
			\clip (-4.3, 0) rectangle (4.3, 3.4);
			\foreach \j in {0,...,5} {
				\draw[thick, ardesia] (-5, \bound + \j*\dplus) -- (5, \bound + \j*\dplus);
			}
			\foreach \i in {-8,...,8} {
				\draw[thick, ardesia] (\i*\dplus, \bound) -- (\i*\dplus, 3.5);
			}
			\foreach \i in {-8,...,8} {
				\foreach \j in {0,...,5} {
					\filldraw[fill=white, draw=black, thick] (\i*\dplus, \bound + \j*\dplus) circle (2pt);
				}
			}
		\end{scope}
		
		\begin{scope}
			\clip (-4.3, 0) rectangle (4.3, -3.4);
			\foreach \j in {0,...,6} {
				\draw[thick, ardesia] (-5, -\bound - \j*\dminus) -- (5, -\bound - \j*\dminus);
			}
			
			\foreach \i in {-10,...,10} {
				\pgfmathsetmacro{\xintB}{\i*\dminus}
				
				\gdef\hasPartner{0}
				\foreach \k in {-8,...,8} {
					\pgfmathsetmacro{\xintT}{\k*\dplus}
					\pgfmathsetmacro{\dist}{abs(\xintT - \xintB)}
					\ifdim \dist cm < 0.25 cm
					\gdef\hasPartner{1}
					\fi
				}
				
				\ifnum\hasPartner=1
				\draw[thick, ardesia] (\xintB, -\bound) -- (\xintB, -3.5);
				\foreach \j in {0,...,6} {
					\filldraw[fill=white, draw=black, thick] (\xintB, -\bound - \j*\dminus) circle (2pt);
				}
				\else
				\draw[thick, ardesia] (\xintB, -0.4) -- (\xintB, -3.5);
				
				\foreach \j in {1,...,6} {
					\filldraw[fill=white, draw=black, thick] (\xintB, -\bound - \j*\dminus) circle (2pt);
				}
				
				\node[text=ottanio, font=\huge, fill=white, inner sep=1pt] at (\xintB, -0.28) {$\boldsymbol{\top}$};
				\fi
			}
		\end{scope}
		
		\draw[very thick, black] (-4.9, 0) node[left] {$\Hh^0_T$} -- (4.7, 0);
		
		\draw[line width=2pt, violet, -{Latex[length=3.5mm, width=2.5mm]}] 
		(1.8, 0.4) -- (3.0, 0.4) node[midway, above=1pt, fill=white, inner sep=1pt] {$\mathbf{b}$};
		
	\end{tikzpicture}
	\caption{Atomistic representation of a semi-coherent interface (van der Merwe misfit model).}
\end{figure}
\subsection{Admissible displacements and elastic energy}
To conclude, the mathematical framework is presented. We fix a defect density function $f \in L^2_T(\R)$ and a constant Burgers vector $\mathbf{b} \in \R^2$. 
\begin{definition}[Displacement field]\label{def: displacement}
Let $U^{+} \in H^1_T(\Hh^+_T,\R^2)$ and $U^{-} \in H^1_T(\Hh^-_T,\R^2)$. We define the global displacement field $U \in SBV^2_T(\R^2,\R^2)$ as
\begin{equation*}
	U \coloneqq U^{+} \mathds{1}_{\Hh^+_T} + U^{-} \mathds{1}_{\Hh^-_T}.
\end{equation*}
\end{definition}
We denote the absolutely continuous part of the gradient as $\beta \coloneqq \nabla U$. Since the dislocation defects are assumed to be distributed along the interface $\Hh^0_T$, the bulk regions $\Hh^+_T \setminus \Hh^0_T$ and $\Hh^-_T \setminus \Hh^0_T$ are free of singularities. This justifies the fact that the strain $\beta$ admits regular Sobolev potentials $U^+$ and $U^-$ within the two half-planes. The presence of defects on the interface results in the jump of the displacement field across $\Hh^0_T$, which translates distributionally into a localized curl constraint. We therefore define the space of admissible displacements associated to a given density $f$ as:
\begin{equation*}
AD(f) \coloneqq \left\{ U \in SBV^2_T(\R^2,\R^2) \suchthat \Curl \beta = f \mathbf{b} \mathcal{H}^1\llcorner\Hh^0_T \text{ in } \mathcal{D}'(\Hh_T) \right\}.
\end{equation*}
Let $G>0$ and $\nu \in (-\frac{1}{2},\frac{1}{2})$ denote, respectively, the shear modulus and the Poisson ratio of the elastic material. The isotropic linear elastic energy of an admissible configuration $U \in AD(f)$ is given by:
\begin{align*}
E_{\mathrm{el}}(U)
&\coloneqq  \frac{1}{2} \int_{\Hh_T} \left(
\frac{2G\nu}{1-2\nu} \, |\operatorname{tr}(e(U))|^2
+ 2G |e(U)|^2 \right) \ud x \ud y \\
&= \frac{1}{2}  \int_{\Hh_T} \sigma(U) : e(U) \, \ud x \ud y ,
\end{align*}
where the symmetric strain tensor is $e(U)\coloneqq \frac{\nabla U + \nabla^T U}{2}$ and the associated Cauchy stress tensor reads
\begin{equation} \label{eq: stress tensore}
\sigma(U) \coloneqq \frac{2G\nu}{1-2\nu} \operatorname{tr}(e(U)) \, \mathrm{Id} + 2G e(U).
\end{equation}

The main objective of our analysis is to characterize the optimal elastic configurations for a prescribed defect distribution. Namely, we aim to study the following minimization problem:
\begin{equation*}
\inf_{U \in AD(f)} E_{\mathrm{el}} (U). \tag{P}
\end{equation*}
\begin{lemma} \label{lemma: euler_lagrange}
Let $U \in AD(f)$ be a minimizer of problem (P) such that its restrictions satisfy $U^\pm \in H^1_T(\Hh_T^\pm, \R^2)$. Then $U$ satisfies the Euler--Lagrange equations:
\begin{equation*}
	\operatorname{div} \sigma(U) = 0 \quad \text{in } \mathcal{D}'(\Hh_T^\pm).
\end{equation*}
\end{lemma}

\begin{proof}
Let $V \in C_c^\infty(\Hh_T^\pm, \R^2)$ and $t \in \R$. Linearity leads us to
\begin{align*}
	E_{\mathrm{el}}(U + tV) &= \frac{1}{2}  \int_{\Hh_T^\pm} \sum_{i,j=1}^{2} \big[ \sigma_{ij}(U) + t\sigma_{ij}(V) \big] \big[ e_{ij}(U) + te_{ij}(V) \big] \ud x \ud y \\
	&= E_{\mathrm{el}}(U) + \frac{t}{2}  \int_{\Hh_T^\pm} \sum_{i,j=1}^{2} \big[ \sigma_{ij}(U) e_{ij}(V) + \sigma_{ij}(V) e_{ij}(U) \big] \ud x \ud y + t^2 E_{\mathrm{el}}(V).
\end{align*}
Evaluating the derivative at $t=0$ and setting it equal to zero, we obtain:
\begin{align*}
	 \left.\frac{\partial}{\partial t} E_{\mathrm{el}}(U + tV)\right|_{t=0} &= \frac{1}{2}  \int_{\Hh_T^\pm} \sum_{i,j=1}^{2} \big[ \sigma_{ij}(U) e_{ij}(V) + \sigma_{ij}(V) e_{ij}(U) \big] \ud x \ud y \\
	&=  \int_{\Hh_T^\pm} \sum_{i,j=1}^{2} \sigma_{ij}(U) e_{ij}(V) \ud x \ud y\\
	&=  \int_{\Hh_T^\pm} \sum_{i,j=1}^{2} \sigma_{ij}(U) \partial_j V_i \ud x \ud y=0\,,
\end{align*}
where in the second equality we have used symmetry of the elasticity tensor ($c_{ijkl} = c_{klij}$) and the identity $\sum_{i,j=1}^{2} \sigma_{ij}(U) e_{ij}(V) = \sum_{i,j=1}^{2} \sigma_{ij}(V) e_{ij}(U)$, while in the third equality we used $\sum_{i,j=1}^{2} \sigma_{ij}(U) e_{ij}(V) = \sum_{i,j=1}^{2} \sigma_{ij}(U) \partial_j V_i$. Now, integrating by parts and using the compact support of \(V\), the result follows by the Fundamental Lemma of the Calculus of Variations.
\end{proof}

\subsection{Elastic Dirichlet-to-Neumann map}
In this subsection, we establish the connection between the 2D bulk elastic energy and the 1D fractional energy on the interface. In particular, we construct the Dirichlet-to-Neumann map associated to the governing equation for the displacement, thereby reducing the 2D problem to a 1D one. We now briefly recall two key conditions that the minimizer of (P) must satisfy.
\begin{hypothesis} \label{hyp: elasticity}
The material is modeled as a linear isotropic elastic continuum in the absence of external body forces. The displacement field $U$ satisfies the equilibrium equations:
\begin{equation} \label{eq: navier_lame}
	\operatorname{div} \sigma(U) = 0 \quad \text{in } \mathcal{D}'(\Hh^\pm_T),
\end{equation}
which corresponds to the Euler--Lagrange equation of \(E_{\mathrm{el}}\), as explained in \Cref{lemma: euler_lagrange}.
Now, rewriting the stress tensor we get
\begin{equation*}
	\sigma(U)= G(\nabla U +\nabla^T U) +\frac{2G\nu}{1-2\nu}(\operatorname{div} U)\operatorname{Id}.
\end{equation*}
Applying \cref{eq: navier_lame}, we obtain
\begin{align}
	0&= G \operatorname{div} (\nabla U) + G \operatorname{div} (\nabla^T U) + \frac{2G\nu}{1-2\nu} \operatorname{div} \left[ \left( \operatorname{div} U \right) \operatorname{Id} \right] \nonumber\\
	&= G \Delta U + G \nabla \left( \operatorname{div} U \right) + \frac{2G\nu}{1-2\nu} \nabla \left( \operatorname{div} U \right) \nonumber \\
	&= G \Delta U + \frac{G}{1-2\nu} \nabla \left( \operatorname{div} U \right) \label{eq: Navier-Lame eqn},
\end{align}
which corresponds to the Navier--Lamé equation of equilibrium in the theory of dislocations (see for example \cite{cinesi20}).
\end{hypothesis}
To properly model an edge dislocation, the minimizer has to satisfy specific parity conditions across the slip plane $\Hh^0_T$, similarly to what has been done in \cite{cinesi20}, depending on the direction of the Burgers vector. This is why it is crucial here to distinguish between the Read--Shockley and van der Merwe models.

\renewcommand{\thehypothesis}{2A}
\renewcommand{\theHhypothesis}{2A}
\begin{hypothesis}[Mirror Symmetry - van der Merwe] \label{hyp: symmetry_VDM}
In particular, in the van der Merwe model, the minimizer satisfies the following symmetries:
\begin{equation*}
	\begin{aligned}
		U_1(x, y) &= -U_1(x, -y), \\
		U_2(x, y) &= U_2(x, -y),
	\end{aligned}
\end{equation*}
for every $(x,y) \in \Hh^+_T$. Consequently, the traces of the displacement on the interface will be denoted as:
\begin{equation*} 
	\begin{aligned}
		U_1(x, 0^+) &\coloneqq \lim_{y \to 0^+} U_1(x,y) = -U_1(x, 0^-) \eqqcolon u(x), \\
		U_2(x, 0^+) &\coloneqq \lim_{y \to 0^+} U_2(x,y) = U_2(x, 0^-) \eqqcolon v(x),
	\end{aligned}
\end{equation*}
where $g(x, 0^\pm) \coloneqq \lim_{y \to 0^\pm} g(x,y)$ denotes the one-sided trace onto the interface $\Hh^0_T$.
\end{hypothesis}

\addtocounter{hypothesis}{-1}
\renewcommand{\thehypothesis}{2B}
\renewcommand{\theHhypothesis}{2B}
\begin{hypothesis}[Mirror Symmetry - Read--Shockley] \label{hyp: symmetry_RS}
On the other hand, in the Read--Shockley model, the minimizer satisfies the following symmetries:
\begin{equation*}
	\begin{aligned}
		U_1(x, y) &= U_1(x, -y), \\
		U_2(x, y) &= -U_2(x, -y),
	\end{aligned}
\end{equation*}
for every $(x,y) \in \Hh^+_T$. Consequently, the traces of the displacement on the interface, in this case, will be denoted as:
\begin{equation*}
	\begin{aligned}
		U_1(x, 0^+) &= U_1(x, 0^-) \eqqcolon v(x), \\
		U_2(x, 0^+) &= -U_2(x, 0^-) \eqqcolon u(x).
	\end{aligned}
\end{equation*}
\end{hypothesis}

Regarding the far-field assumptions for both models, since the problem is set on the periodic strip \(\Hh_T\), we assume the slip profile $u \in H^{\frac{1}{2}}_T(\R)$ to have zero average, i.e., $\fint_0^T u(x) \ud x = 0$, ensuring the finiteness of the $L^2$-norm of the displacement gradient in the unbounded half-spaces (see \cref{rmk: media nulla}). 

We point out that \cref{hyp: elasticity} and \cref{hyp: symmetry_VDM} (or \cref{hyp: symmetry_RS}) need not be assumed \textit{a priori}: by \Cref{lemma: euler_lagrange} and the convexity argument in \cref{rmk: simmetria del minimo}, these conditions are naturally satisfied by the minimizers.
\renewcommand{\thehypothesis}{\arabic{hypothesis}}
\renewcommand{\theHhypothesis}{\arabic{hypothesis}}
\begin{remark}\label{rmk: riduzione a H^+}
Without loss of generality, we can restrict our minimization problem to symmetric displacement fields satisfying \cref{hyp: symmetry_VDM} or \cref{hyp: symmetry_RS}. For such configurations, the energy is symmetric with respect to the interface \(\Hh^0_T\). Consequently, the total elastic energy over the periodic strip can be rewritten as:
\begin{align*}
	E_{\mathrm{el}}(U) &=  \int_{\Hh^+_T} \left( \frac{2G\nu}{1-2\nu} \, |\operatorname{tr}(e(U))|^2 + 2G |e(U)|^2 \right) \ud x \ud y \\
	&=  \int_{\Hh^+_T} \sigma(U) : e(U) \, \ud x \ud y.
\end{align*}   
\end{remark}
The next result will be formulated in the van der Merwe setting (\cref{hyp: symmetry_VDM}).
\begin{definition}[Elastic extension] \label{def: elastic_extension}
Given $u \in H^{\frac{1}{2}}_{T}(\R)$ with zero average, we define its \textit{elastic extension} to the periodic strip $\Hh_T$ as the unique displacement field $U \in SBV^2_T(\R^2,\R^2)$ which solves
\begin{equation*}
	\begin{cases}
		\operatorname{div} \sigma(U) = 0 \qquad &\text{in } \Hh^\pm_T\,, \\
		U_1(x,y) = -U_1(x,-y) \qquad &\text{for } y > 0\,, \\
		U_2(x,y) = U_2(x,-y) \qquad &\text{for } y > 0\,, \\
		U_1(x,0^+) = u(x) \qquad &\text{on } \Hh^0_T\,,\\
		\sigma_{22}(x,0^+) = \sigma_{22}(x,0^-) \qquad &\text{on } \Hh^0_T\,,
	\end{cases}
\end{equation*}
where \(\sigma\) is the stress tensor defined in \cref{eq: stress tensore}. The interface condition corresponds to the natural continuity of the normal traction across the interface (cf.~\cite{cinesi20}).
\end{definition}

\begin{theorem}[Elastic Dirichlet-to-Neumann map] \label{thm: elastic_dtn}
Let $u \in H^{\frac{1}{2}}_T(\R)$ with zero average and let $U \in SBV^2_T(\R^2,\R^2)$ be its elastic extension on \(\Hh_T\). Then, the normal stress satisfies
\begin{equation} \label{eq: zero_normal_stress}
	\sigma_{22}(x, 0^+) = \sigma_{22}(x, 0^-) = 0 \quad \text{for every } x \in [0,T],
\end{equation}
and the following identity holds:
\begin{equation} \label{eq: energy_reduction_elastic}
	\frac{1}{2} \int_{\Hh^+_T} \sigma(U) : e(U) \ud x \ud y = \frac{G}{1-\nu} \frac{1}{4\pi} [u]_{\frac{1}{2}, T}^2.
\end{equation}
\end{theorem}
\begin{remark}
Before giving the proof of the result, we recall here for the reader's convenience the components of the stress operator $\sigma(U)$ on the two-dimensional displacement field $U(x,y) = (U_1(x,y), U_2(x,y))^T$.
\begin{align}
	\sigma_{11}(U) &= 2G \partial_x U_1 + \frac{2G\nu}{1-2\nu} (\partial_x U_1 + \partial_y U_2) = \frac{G}{1-2\nu} \left[ 2(1-\nu)\partial_x U_1 + 2\nu\partial_y U_2 \right], \nonumber \\
	\sigma_{22}(U) &= 2G \partial_y U_2 + \frac{2G\nu}{1-2\nu} (\partial_x U_1 + \partial_y U_2) = \frac{G}{1-2\nu} \left[ 2(1-\nu)\partial_y U_2 + 2\nu\partial_x U_1 \right], \label{eq: sigma_22}\\
	\sigma_{12}(U) &= 2G \left( \frac{\partial_y U_1 + \partial_x U_2}{2} \right) = G(\partial_y U_1 + \partial_x U_2) . \nonumber
\end{align}
\end{remark}
\begin{proof}
The proof adopts a Fourier approach.
\subsubsection*{Step 1: Fourier series expansion and substitution into the Navier--Lamé system}
Let $\omega_n = \frac{2\pi n}{T}$ be the angular frequency associated to the $n$-th Fourier mode for $n \in \mathbb{Z}$. From the classical theory of Fourier series, the expansion of each component of the displacement field $U(x,y) = (U_1(x,y), U_2(x,y))^T$ in the upper half-plane $\Hh^+_T$ (namely, for $y > 0$) reads as:
\begin{equation*}
	U_1(x,y) = \sum_{n \in \mathbb{Z}} \hat{U}_{1,n}(y) e^{i\omega_n x}, \quad \text{and} \quad U_2(x,y) = \sum_{n \in \mathbb{Z}} \hat{U}_{2,n}(y) e^{i\omega_n x},
\end{equation*}
where the coefficients $\hat{U}_{1,n}(y)$ and $\hat{U}_{2,n}(y)$ are real- or complex-valued functions. The partial derivatives of each component, which appear in \cref{eq: Navier-Lame eqn}, take the form
\begin{align*}
	\partial_x U_1(x,y) &= \sum_{n \in \mathbb{Z}} i\omega_n \hat{U}_{1,n}(y) e^{i\omega_n x}, \\
	\partial^2_{x} U_1(x,y) &= \sum_{n \in \mathbb{Z}} (i\omega_n)^2 \hat{U}_{1,n}(y) e^{i\omega_n x} = \sum_{n \in \mathbb{Z}} -\omega_n^2 \hat{U}_{1,n}(y) e^{i\omega_n x}, \\
	\partial^2_{y} U_1(x,y) &= \sum_{n \in \mathbb{Z}} \partial^2_{y} \hat{U}_{1,n}(y) e^{i\omega_n x}, \\
	\partial^2_{xy} U_2(x,y) &= \partial_x \left( \sum_{n \in \mathbb{Z}} \partial_y \hat{U}_{2,n}(y) e^{i\omega_n x} \right) = \sum_{n \in \mathbb{Z}} i\omega_n \partial_y \hat{U}_{2,n}(y) e^{i\omega_n x}.
\end{align*}
Analogously, for the second component of the displacement field $U_2(x,y)$, we obtain
\begin{align*}
	\partial^2_{x} U_2(x,y) &= \sum_{n \in \mathbb{Z}} -\omega_n^2 \hat{U}_{2,n}(y) e^{i\omega_n x}, \\
	\partial^2_{y} U_2(x,y) &= \sum_{n \in \mathbb{Z}} \partial^2_{y} \hat{U}_{2,n}(y) e^{i\omega_n x}, \\
	\partial^2_{xy} U_1(x,y) &= \sum_{n \in \mathbb{Z}} i\omega_n \partial_y \hat{U}_{1,n}(y) e^{i\omega_n x}.
\end{align*}
Given these expressions, we proceed to plug them into \cref{eq: navier_lame}, obtaining the following system of equations
\begin{equation*}
	\begin{cases}
		2(1-\nu)\partial^2_{x} U_1 + (1-2\nu)\partial^2_{y} U_1 + \partial^2_{xy} U_2 = 0, \\
		(1-2\nu)\partial^2_{x} U_2 + 2(1-\nu)\partial^2_{y} U_2 + \partial^2_{xy} U_1 = 0,
	\end{cases}
\end{equation*}
which, substituting the respective Fourier series, takes the form:
\begin{equation*}
	\begin{cases}
		\sum_{n \in \mathbb{Z}} \left[ -2(1-\nu)\omega_n^2 \hat{U}_{1,n}(y) + (1-2\nu)\partial^2_{y} \hat{U}_{1,n}(y) + i\omega_n \partial_y \hat{U}_{2,n}(y) \right] e^{i\omega_n x} = 0, 			\vspace{3pt} \\
		\sum_{n \in \mathbb{Z}} \left[ -(1-2\nu)\omega_n^2 \hat{U}_{2,n}(y) + 2(1-\nu)\partial^2_{y} \hat{U}_{2,n}(y) + i\omega_n \partial_y \hat{U}_{1,n}(y) \right] e^{i\omega_n x} = 0.
	\end{cases}
\end{equation*}
By the orthogonality of the Fourier basis $\{e^{i\omega_n x}\}_{n \in \mathbb{Z}}$ in the space $L^2(0,T)$, each coefficient of the series must vanish for every $n \in \mathbb{Z}$. We thus obtain the following coupled system of ODEs with respect to the variable $y$:
\begin{align}
	(1-2\nu)\partial^2_{y} \hat{U}_{1,n}(y) - 2(1-\nu)\omega_n^2 \hat{U}_{1,n}(y) + i\omega_n \partial_y \hat{U}_{2,n}(y) &= 0, \label{eq: ode_accoppiata1} \\
	2(1-\nu)\partial^2_{y} \hat{U}_{2,n}(y) - (1-2\nu)\omega_n^2 \hat{U}_{2,n}(y) + i\omega_n \partial_y \hat{U}_{1,n}(y) &= 0. \label{eq: ode_accoppiata2}
\end{align}
\subsubsection*{Step 2: Decoupling of the system and derivation of the fourth-order equation}
In order to solve the system, we isolate the coupling term containing $\hat{U}_{2,n}$ from \cref{eq: ode_accoppiata1}:
\begin{equation} \label{eq: isolamento_u2}
	i\omega_n \partial_y \hat{U}_{2,n}(y) = 2(1-\nu)\omega_n^2 \hat{U}_{1,n}(y) - (1-2\nu)\partial^2_{y} \hat{U}_{1,n}(y).
\end{equation}
Taking the second ordinary differential \cref{eq: ode_accoppiata2} and differentiating each term with respect to $y$ yields:
\begin{align*}
	0&=\partial_y \left[ 2(1-\nu)\partial^2_{y} \hat{U}_{2,n}(y) - (1-2\nu)\omega_n^2 \hat{U}_{2,n}(y) + i\omega_n \partial_y \hat{U}_{1,n}(y) \right]\\
	&=2(1-\nu)\partial^3_{y} \hat{U}_{2,n}(y) - (1-2\nu)\omega_n^2 \partial_y \hat{U}_{2,n}(y) + i\omega_n \partial^2_{y} \hat{U}_{1,n}(y).
\end{align*}
We then multiply the entire equation by the complex constant $i\omega_n$ (recalling that $(i\omega_n)^2 = -\omega_n^2$)
\begin{equation} \label{eq: ode2_moltiplicata}
	2(1-\nu)\partial^2_{y} \left( i\omega_n \partial_y \hat{U}_{2,n}(y) \right) - (1-2\nu)\omega_n^2 \left( i\omega_n \partial_y \hat{U}_{2,n}(y) \right) - \omega_n^2 \partial^2_{y} \hat{U}_{1,n}(y) = 0,
\end{equation}
and substituting \cref{eq: isolamento_u2} into \cref{eq: ode2_moltiplicata} gives:
\begin{align} 
	2(1-\nu)\partial^2_{y} \left[ 2(1-\nu)\omega_n^2 \hat{U}_{1,n}(y) - (1-2\nu)\partial^2_{y} \hat{U}_{1,n}(y) \right] \nonumber \\
	- (1-2\nu)\omega_n^2 \left[ 2(1-\nu)\omega_n^2 \hat{U}_{1,n}(y) - (1-2\nu)\partial^2_{y} \hat{U}_{1,n}(y) \right] - \omega_n^2 \partial^2_{y} \hat{U}_{1,n}(y) = 0.\label{eq: ode_combined}
\end{align}
Now, multiplying and collecting the terms according to the order of differentiation of the function $\hat{U}_{1,n}(y)$ and simplifying the coefficients, we obtain
\begin{align*}
	\left[ -2(1-\nu)(1-2\nu) \right] \partial^4_{y} \hat{U}_{1,n}(y)& \nonumber \\
	- 2\omega_n^2 \left[ -2(1-\nu)(1-2\nu) \right] \partial^2_{y} \hat{U}_{1,n}(y)& \nonumber \\
	+ \omega_n^4 \left[ -2(1-\nu)(1-2\nu) \right] \hat{U}_{1,n}(y) &= 0.
\end{align*}
Since the Poisson ratio satisfies $\nu \in (-\frac{1}{2}, \frac{1}{2})$, the quantity $-2(1-\nu)(1-2\nu)$ is strictly non-zero: dividing then \cref{eq: ode_combined} by this constant, we arrive at the decoupled fourth-order ordinary differential equation:
\begin{equation} \label{eq: ode_quarto_ordine}
	\partial^4_{y} \hat{U}_{1,n}(y) - 2\omega_n^2 \partial^2_{y} \hat{U}_{1,n}(y) + \omega_n^4 \hat{U}_{1,n}(y) = 0.
\end{equation}
\subsubsection*{Step 3: Resolution of the ODE and energy finiteness constraint in the half-plane}
The characteristic equation associated with \cref{eq: ode_quarto_ordine} admits two symmetric real roots, each having an algebraic multiplicity equal to two: $\lambda_1 = \omega_n$ and $\lambda_2 = -\omega_n$. To simultaneously treat the cases where $n > 0$ and $n < 0$, we express the roots in terms of the modulus of the frequency, namely $\pm |\omega_n|$.
The general solution takes the form:
\begin{equation*}
	\hat{U}_{1,n}^+(y) = (A_n^+ + B_n^+ y)e^{|\omega_n|y} + 	(C_n^+ + D_n^+ y)e^{-|\omega_n|y}.
\end{equation*}
To ensure finiteness of the $L^2$-norm of the gradient, we impose \(	A_n^+ = 0\) and \( B_n^+ = 0\) for every \(n \in \Z \setminus\{0\}\).
Relabeling the free constants for simplicity, the admissible configuration for the horizontal and vertical displacements in the upper half-plane $y > 0$ is written as:
\begin{align}
	\hat{U}_{1,n}^+(y) &= (A_n^+ - B_n^+ |\omega_n| y) e^{-|\omega_n|y}, \label{eq: sol_generale1} \\
	\hat{U}_{2,n}^+(y) &= -i \operatorname{sgn}(\omega_n) (C_n^+ - D_n^+ |\omega_n| y) e^{-|\omega_n|y}, \label{eq: sol_generale2}
\end{align}
where $A_n^+, B_n^+, C_n^+, D_n^+ \in \mathbb{C}$ for every $n \in \mathbb{Z} \setminus \{0\}$.
\subsubsection*{Step 4: Internal coupling and determination of the coefficients}
At this stage, however, for \cref{eq: sol_generale1,eq: sol_generale2} to actually be solutions of the original Navier--Lamé system, they must satisfy \cref{eq: ode_accoppiata1,eq: ode_accoppiata2}. We compute the first and second derivatives with respect to $y$ of the solution $\hat{U}_{1,n}^+(y)$:
\begin{align}
	\partial_y \hat{U}_{1,n}^+(y) &= -|\omega_n| \left[ (A_n^+ + B_n^+) - B_n^+ |\omega_n| y \right] e^{-|\omega_n|y}, \nonumber \\
	\partial^2_{y} \hat{U}_{1,n}^+(y) &= \omega_n^2 \left[ (A_n^+ + 2B_n^+) - B_n^+ |\omega_n| y \right] e^{-|\omega_n|y}. \label{eq: derivata2_yU_1_n}
\end{align}
Analogously, for $\hat{U}_{2,n}^+(y)$:
\begin{align*}
	\partial_y \hat{U}_{2,n}^+(y) &= i \operatorname{sgn}(\omega_n) |\omega_n| \left[ (C_n^+ + D_n^+) - D_n^+ |\omega_n| y \right] e^{-|\omega_n|y}. 
\end{align*}
In the Navier--Lamé system, this term appears multiplied by $i\omega_n$. Since $\omega_n \operatorname{sgn}(\omega_n) = |\omega_n|$, we observe that $i\omega_n \cdot i \operatorname{sgn}(\omega_n) |\omega_n| = i^2 |\omega_n|^2 = -\omega_n^2$. Hence:
\begin{equation}
	i\omega_n \partial_y \hat{U}_{2,n}^+(y) = -\omega_n^2 \left[ (C_n^+ + D_n^+) - D_n^+ |\omega_n| y \right] e^{-|\omega_n|y}. \label{eq: iomegaderivata_yU_2_n}
\end{equation}
Substituting \cref{eq: sol_generale1,eq: derivata2_yU_1_n,eq: iomegaderivata_yU_2_n} into \cref{eq: ode_accoppiata1} and simplifying the expression, we obtain
\begin{equation*}
	-2(1-\nu) \left[ A_n^+ - B_n^+ |\omega_n| y \right] + (1-2\nu) \left[ (A_n^+ + 2B_n^+) - B_n^+ |\omega_n| y \right] - \left[ (C_n^+ + D_n^+) - D_n^+ |\omega_n| y \right] = 0.
\end{equation*}
Since the equation must hold for every $y>0$, the coefficients of the linear term in $|\omega_n|y$ and of the constant term must vanish separately. 
Equating the coefficients to zero, we obtain
\begin{equation}\label{eq: D_n relazione B_n}
	D_n^+ = -B_n^+,
\end{equation}
and 
\begin{equation}\label{eq: C_n relazione A_n e B_n}
	C_n^+=-A_n^+ + (3-4\nu)B_n^+.
\end{equation}
\subsubsection*{Step 5: Symmetry conditions and vanishing of the normal stress}
We now analyze the behavior at the interface $y = 0$; we refer to \cref{rmk: leggitimita_pass_elastico} for a justification of the subsequent steps. Using \cref{hyp: symmetry_VDM}, this immediately implies that:
\begin{equation*}
	A_n^+ = -A_n^- \quad \text{and} \quad C_n^+ = C_n^-.
\end{equation*}
By \cref{eq: sigma_22}, we obtain
\begin{align}
	\hat{\sigma}_{22,n}^+(0) &= \frac{G}{1-2\nu} \left[ 2(1-\nu)\partial_y \hat{U}_{2,n}^+(0) + 2\nu(i\omega_n) \hat{U}_{1,n}^+(0) \right] \nonumber\\
	&= \frac{i\omega_n G}{1-2\nu} \left[ 2(1-\nu)(C_n^+ + D_n^+) + 2\nu A_n^+ \right] \nonumber\\
	&= i\omega_n G \left[ -2A_n^+ + 4(1-\nu)B_n^+ \right] \label{eq: sigma_22_n},
\end{align}
where in the second equality we have used that $\hat{U}_{1,n}^+(0) = A_n^+$ and that $\partial_y \hat{U}_{2,n}^+(0) = i\omega_n (C_n^+ + D_n^+)$, and in the last one we have used \cref{eq: D_n relazione B_n,eq: C_n relazione A_n e B_n}.	
By the interface transmission condition in \cref{def: elastic_extension}, the normal stress is continuous across the interface, which in Fourier components reads $\hat{\sigma}_{22,n}^+(0) = \hat{\sigma}_{22,n}^-(0)$. Furthermore, from \cref{hyp: symmetry_VDM} on $U_1$ and $U_2$, the normal stress exhibits an odd symmetry with respect to the interface, namely $\hat{\sigma}_{22,n}^+(0) = -\hat{\sigma}_{22,n}^-(0)$. The only solution compatible with both conditions is the trivial one:
\begin{equation*}
	\hat{\sigma}_{22,n}^+(0) = 0.
\end{equation*}
This proves \cref{eq: zero_normal_stress}. Hence, from \cref{eq: sigma_22_n}, we obtain
\begin{equation} \label{eq: A_n relazione B_n}
	B_n^+ = \frac{A_n^+}{2(1-\nu)}.
\end{equation}
\subsubsection*{Step 6: Computation of the shear stress and reduction of the elastic energy}
The shear stress evaluated at the interface $y=0$ is:
\begin{align*}
	\hat{\sigma}_{12,n}^+(0) &= G \left( \partial_y \hat{U}_{1,n}^+(0) + i\omega_n \hat{U}_{2,n}^+(0) \right) \nonumber\\
	&= G |\omega_n| \left( -A_n^+ - B_n^+ + C_n^+ \right)\\
	&= G |\omega_n| \left( -2A_n^+ + 2(1-2\nu)B_n^+ \right)\\
	&= G |\omega_n| A_n^+ \left( \frac{-1}{1-\nu} \right),
\end{align*}
where we have used \cref{eq: C_n relazione A_n e B_n,eq: A_n relazione B_n}.	
Finally, remembering that \(U_1(x,0^+)=u(x)\), we have that \(A_n^+=\hat u_n\) for every \(n \in \Z \setminus \{0\}\). Thus,
\begin{equation*}
	\hat{\sigma}_{12,n}^+(0) = - \frac{G}{1-\nu} |\omega_n| \hat{u}_n.
\end{equation*}
We are now ready to compute the elastic energy in the upper half-plane. Integration by parts yields
\begin{equation*}
	\int_{\Hh^+_T} \sigma(U) : e(U) \ud x \ud y = \int_{\partial \Hh^+_T} \left( \sigma \cdot \mathbf{n}^+ \right) \cdot U \ud \sigma.
\end{equation*}
Now, arguing similarly as in \cref{eq: identita energia}, we get:
\begin{align*}
	\frac{1}{2}  \int_{\Hh^+_T} \sigma(U) : e(U) \ud x \ud y &= - \frac{1}{2}  \int_0^T \left( \sigma_{12}(x,0^+) U_1(x,0^+) + \sigma_{22}(x,0^+) U_2(x,0^+) \right) \ud x\\
	&=-\frac{1}{2}  \int_0^T \sigma_{12}(x,0^+) u(x) \ud x\\
	&= -\frac{T}{2}  \sum_{n \in \mathbb{Z} \setminus \{0\}} \hat{\sigma}_{12,n}^+(0) \overline{\hat{u}_n}\\
	&=-\frac{T}{2} \sum_{n \in \mathbb{Z} \setminus \{0\}} \left( - \frac{G}{1-\nu} |\omega_n| \hat{u}_n \right) \overline{\hat{u}_n} \\
	&= \frac{G}{1-\nu} \frac{T}{2} \sum_{n \in \mathbb{Z} \setminus \{0\}} |\omega_n| |\hat{u}_n|^2\\
	&= \frac{G}{1-\nu} \frac{1}{4\pi} [u]_{\frac{1}{2}, T}^2,
\end{align*}
where in the second equality we have used \cref{eq: zero_normal_stress}, in the third one Parseval identity and for the last one see \cref{rmk: h12infourier}. This proves \cref{eq: energy_reduction_elastic}.
\end{proof}

\begin{remark}
Regarding the Read--Shockley setting, having required \cref{hyp: symmetry_RS}, we will have that \(\sigma_{12}(x,0^+)=0\) for a.e. \(x \in \R\), so that integration by parts will lead to
\begin{align*}
	\frac{1}{2}\int_{\Hh^+_T} \sigma(U) : e(U) \ud x \ud y &= - \frac{1}{2} \int_0^T \left( \sigma_{12}(x,0^+) U_1(x,0^+) + \sigma_{22}(x,0^+) U_2(x,0^+) \right) \ud x\\
	&=- \frac{1}{2} \int_0^T \sigma_{22}(x,0^+) u(x) \ud x\\
	&= -\frac{T}{2} \sum_{n \in \mathbb{Z} \setminus \{0\}} \hat{\sigma}_{22,n}^+(0) \overline{\hat{u}_n}\\
	&=-\frac{T}{2} \sum_{n \in \mathbb{Z} \setminus \{0\}} \left( - \frac{G}{1-\nu} |\omega_n| \hat{u}_n \right) \overline{\hat{u}_n} \\
	&= \frac{G}{1-\nu} \frac{T}{2} \sum_{n \in \mathbb{Z} \setminus \{0\}} |\omega_n| |\hat{u}_n|^2\\
	&= \frac{G}{1-\nu} \frac{1}{4\pi} [u]_{\frac{1}{2}, T}^2,
\end{align*}
yielding the same boundary energy. 
\end{remark}

\begin{remark} \label{rmk: leggitimita_pass_elastico}
	Let \(\mathbf{n}=(0,-1)\). The evaluation of the stress components at the interface $y=0$ is justified by functional analysis arguments, analogously to the scalar case discussed in \Cref{lemma: normal_trace_harmonic}. For any displacement field $U \in H^1_T(\Hh^+_T, \R^2)$ satisfying $\operatorname{div} \sigma(U) = 0$ in the weak sense, the normal trace $\sigma(U)\mathbf{n} = (-\sigma_{12}(\cdot, 0^+), -\sigma_{22}(\cdot, 0^+))^T$ can be uniquely extended to a continuous linear operator from the space of elastic extensions in $H^1_T(\Hh^+_T, \R^2)$ (see \cref{def: elastic_extension}) to the dual space $H^{-\frac{1}{2}}_T(\Hh^0_T, \R^2)$. It is defined by duality in the following way
	\begin{equation*}
		\langle \sigma(U)\mathbf{n}, \Psi \rangle_{H^{-\frac{1}{2}}_T \times H^{\frac{1}{2}}_T} \coloneqq \int_{\Hh^+_T} \sigma(U) : e(\Phi) \ud x \ud y,
	\end{equation*}
	for every boundary datum $\Psi \in H^{\frac{1}{2}}_T(\Hh^0_T, \R^2)$ and for every extension $\Phi \in H^1_{\Psi,T}(\Hh^+_T, \R^2)$. The continuity of this operator with respect to the $H^1$ topology guarantees that the normal trace and the infinite Fourier sums can be interchanged.
\end{remark}

\section{Minimization problem} \label{sec: section3} \label{sec: min}
In this section we prove the existence (and uniqueness) of the minimizer for problem (P). The following lemma simplifies the proof of the main theorem. In particular, it explains how the curl constraint can be translated into a condition on the traces. This will imply an admissibility condition for the curl.

\begin{lemma}\label{lemma: trace}
Let \(f_1,f_2 \in L^2_{T}(\R)\), let \(U \in SBV^2_T(\R^2, \R^2)\) be a displacement field such that its restrictions \(U^\pm \in H^1_T(\Hh^\pm_T, \R^2)\), and let \(u^\pm \in H^{\frac{1}{2}}_T(\Hh^0_T, \R^2)\) be, respectively, their traces on the interface \(\Hh^0_T\). Then \(\Curl \beta =(f_1,f_2)\,\mathcal H^1\llcorner\Hh_T^0\) if and only if \( w\coloneqq u^+-u^- \in H^1_T(\Hh^0_T, \R^2) \) and \(-w'= (f_1,f_2)\) in the distributional sense. In particular, for any $f \in L^2_{T}(\R)$, it holds that $U \in AD(f)$ if and only if \(-w'= f \mathbf{b}\).
\end{lemma}

\begin{proof}
We prove the result by density. Let us first assume that \( U^\pm\in C^\infty(\overline{\Hh_T^\pm},\R^2)\), let \(U\) be the displacement field as in \cref{def: displacement}, and let
\(\varphi\in C_c^\infty(\Hh_T,\R^2)\).
Let us recall that 
\begin{equation*}
	\beta \coloneqq  \nabla U^{+} \mathds{1}_{\Hh^+_T} + \nabla U^{-} \mathds{1}_{\Hh^-_T}.
\end{equation*}
Therefore, we obtain
\begin{align*}
	\langle \Curl\beta,\varphi\rangle
	&=
	\int_{\Hh_T^+}
	\bigl(
	\partial_x U^+\cdot\partial_y\varphi
	-
	\partial_y U^+\cdot\partial_x\varphi
	\bigr)\,\ud x\,\ud y
	\\
	&\quad+
	\int_{\Hh_T^-}
	\bigl(
	\partial_x U^-\cdot\partial_y\varphi
	-
	\partial_y U^-\cdot\partial_x\varphi
	\bigr)\,\ud x\,\ud y
	\\
	&\eqqcolon \mathcal I^+ + \mathcal I^- .
\end{align*}
To evaluate \(\mathcal I^+\), we integrate by parts. Denoting by
\(\mathbf{n}^+=(n_x^+,n_y^+)=(0,-1)\) the outer normal to
\(\Hh_T^+\) along the interface, we obtain
\begin{align*}
	\mathcal I^+
	&=
	\int_{\Hh_T^+}
	\bigl(
	\partial_x U^+\cdot\partial_y\varphi
	-
	\partial_y U^+\cdot\partial_x\varphi
	\bigr)\,\ud x\,\ud y
	\\
	&=
	\int_{\Hh_T^+}
	U^+\cdot
	\bigl(
	\partial_y\partial_x\varphi
	-
	\partial_x\partial_y\varphi
	\bigr)\,\ud x\,\ud y +
	\int_{\Hh_T^0}
	U^+\cdot
	\bigl(
	\partial_y\varphi\,n_x^+
	-
	\partial_x\varphi\,n_y^+
	\bigr)\,\ud x
	\\
	&=
	\int_{\Hh_T^0}u^+\cdot\partial_x\varphi\,\ud x ,
\end{align*}
where \(u^+\) denotes the trace of \(U^+\) on \(\Hh_T^0\). Hence,
\begin{equation}\label{eq: H+ part}
	\mathcal I^+
	=
	\int_{\Hh_T^0}u^+\cdot\partial_x\varphi\,\ud x .
\end{equation}
Analogously, on \(\Hh_T^-\), using the outer normal
\(\mathbf{n}^-=(n_x^-,n_y^-)=(0,1)\), we obtain
\begin{equation}\label{eq: H- part}
	\mathcal I^-
	=
	-\int_{\Hh_T^0}u^-\cdot\partial_x\varphi\,\ud x .
\end{equation}
Combining \cref{eq: H+ part,eq: H- part} and setting \(w \coloneqq u^+ - u^-\), we obtain
\begin{equation*}
	\langle \Curl\beta,\varphi\rangle_{\mathcal D'\times\mathcal D}
	=
	\int_{\Hh_T^0}
	w(x) \cdot\partial_x\varphi(x,0)\,\ud x ,
\end{equation*}
which shows that \(\Curl \beta = -w' \mathcal{H}^1\llcorner\Hh_T^0\) in the distributional sense.
As a consequence, the condition  \(\Curl \beta = (f_1,f_2)\,\mathcal H^1\llcorner\Hh_T^0\), for $f_1,f_2\in L^{2}_T(\R)$ holds if and only if
\[
	-w' = (f_1,f_2)
\]
in the distributional sense. This concludes the proof of the equivalence in the smooth case.
We now conclude by density. Let
\( U_n^\pm \in C^\infty(\overline{\Hh_T^\pm},\R^2) \) be smooth
approximations such that
\[
	U_n^\pm\to U^\pm
	\qquad\text{strongly in } H^1_T(\Hh_T^\pm,\R^2),
\]
obtained, for instance, by reflection and mollification. Let \(u_n^\pm\)
be their traces, \(w_n\coloneqq u_n^+-u_n^-\), let \((f_{1,n},f_{2,n})\) denote the
corresponding interfacial datum, and let $\beta_n$ be the corresponding field.  By the smooth case,
\[
	-w_n'= (f_{1,n},f_{2,n}).
\]
Moreover, by continuity of the trace operator,
\(w_n\to w\) in \(H_T^{\frac12}(\Hh_T^0,\R^2)\), while
\(\beta_n\to \beta \) in \(L^2\). Since
\[
	\Curl:L^2(\Hh_T,\R^{2\times2})\longrightarrow
	H^{-1}_T(\Hh_T,\R^2)
\]
is continuous, we have
\[
	(f_{1,n},f_{2,n})\,\mathcal H^1\llcorner\Hh_T^0
	=
	\Curl\beta_n
	\longrightarrow
	\Curl\beta
	=
	(f_1,f_2)\,\mathcal H^1\llcorner\Hh_T^0
\]
strongly in \(H^{-1}_T(\Hh_T,\R^2)\). Passing to the limit in
\(-w_n'=(f_{1,n},f_{2,n}) \) yields
\[
	-w'=(f_{1},f_2)
\]
in the distributional sense. The converse follows similarly, since the approximation argument relies entirely on a chain of equalities.
\end{proof}

\begin{remark}\label{rmk: AD convexity}
The space \(AD(f)\) is convex. Indeed, \(SBV^2_T(\R^2,\R^2)\) is a vector space and, for every \(U,V\in AD(f)\) and \(\lambda\in[0,1]\),
\[
\lambda U+(1-\lambda)V\in SBV^2_T(\R^2,\R^2).
\]
Moreover, by linearity of the \(\Curl\) operator,
\[
\Curl\bigl(\lambda\nabla U+(1-\lambda)\nabla V\bigr)
=f\mathbf b\,\mathcal H^1\llcorner\Hh_T^0.
\]
Hence \(\lambda U+(1-\lambda)V\in AD(f)\).
\end{remark}

\begin{remark} \label{rmk: simmetria del minimo}
Given \(U \in AD(f)\), one can always construct a symmetric competitor \(\Tilde{U} \in AD(f)\) satisfying \cref{hyp: symmetry_VDM} (the same holds for \cref{hyp: symmetry_RS}).
Supposing without loss of generality \(\mathbf{b} = \mathbf{e_1}\), we can define the competitor displacement field \(U^*\) as:
\begin{equation*}
	U^*(x,y) \coloneqq \begin{pmatrix} -U_1(x,-y) \\ U_2(x,-y) \end{pmatrix}.
\end{equation*}
By parity, the vertical jump across \(\Hh^0_T\) is zero and the jump of the reflected field becomes \([U^*_1] = -U_1(x,0^-) - (-U_1(x,0^+)) = [U_1]\). Thus \(U^* \in AD(f)\).	We can then define the symmetrized field:
\begin{equation*}
	\Tilde{U} \coloneqq \frac{U+U^*}{2}.
\end{equation*}
By construction, \(\Tilde{U}\) strictly satisfies the symmetry conditions formulated in \cref{hyp: symmetry_VDM}. The convexity of the space \(AD(f)\) (see \cref{rmk: AD convexity}) guarantees that \(\Tilde{U}\) is still admissible. Furthermore, since the elastic strain energy density is invariant under this coordinate reflection, we have \(E_{\mathrm{el}}(U) = E_{\mathrm{el}}(U^*)\). Finally, the convexity of the elastic energy \(E_{\mathrm{el}}\) immediately yields \(E_{\mathrm{el}}(\Tilde{U}) \leq E_{\mathrm{el}}(U)\).
\end{remark}

\begin{remark}
We point out that for any displacement field $U \in SBV^2_T(\R^2, \R^2)$ such that its restrictions satisfy $U^\pm \in H^1_T(\Hh^\pm_T, \R^2)$, the total elastic energy $E_{\mathrm{el}}(U)$ is finite. Indeed, by standard pointwise inequalities, $|e(U)|^2 \leq |\nabla U|^2$ and $|\operatorname{tr}(e(U))|^2 \leq 2 |\nabla U|^2$. Therefore, there exists a constant $C \coloneqq C(G,\nu)>0$ such that $E_{\mathrm{el}}(U) \leq C \|\nabla U\|_{L^2(\Hh_T)}^2$.
\end{remark}

\begin{theorem} \label{thm: energia elastica e' la seminorma}
Let \(f \in L^2_T(\R)\) with \(\fint f = 0\), and let \(\mathbf{b} \in \{\mathbf{e_1},\mathbf{e_2}\}\). Then, the minimization problem (P)
admits a unique minimizer, i.e., there exists a unique \(\tilde{U} \in AD(f)\) such that 
\[
E_{\mathrm{el}}(\tilde{U}) = \inf_{\substack{U \in AD(f)}} E_{\mathrm{el}}(U).
\]
Moreover, the following holds
\begin{equation} \label{eq: l'en minima e' la seminorma}
	E_{\mathrm{el}}(\tilde{U}) = c(G,\nu) [u]^2_{\frac{1}{2},T}\,,
\end{equation}
where \(c(G,\nu) = \frac{G}{2\pi(1-\nu)}\) and \(u\) is given by
\begin{equation} \label{eq: primitive_u}
	u(x)= -\frac{1}{2} \int^x_0 f(t) \ud t + \frac{1}{2T} \int_0^T \int_0^s f(t) \ud t \ud s\,.		
\end{equation}
\end{theorem}
\begin{proof}
Uniqueness follows from the strict convexity of the elastic energy, together with the prescribed normalization. We prove the case \(\mathbf{b}=\mathbf{e_1}\), while the case \(\mathbf{b}=\mathbf{e_2}\) follows analogously by using \cref{hyp: symmetry_RS}. In view of \cref{rmk: simmetria del minimo}, the minimization problem can be restricted to displacement fields satisfying the symmetries in \cref{hyp: symmetry_VDM}. For such configurations, the reduction described in \cref{rmk: riduzione a H^+} applies.

By \Cref{lemma: trace}, any admissible displacement must satisfy the jump condition \(-(U_1(x, 0^+) - U_1(x, 0^-))' = f(x)\). Using the mirror symmetry, we have \(U_1(x, 0^+) = -U_1(x, 0^-)\), which implies that the trace on the upper boundary, \(u(x) \coloneqq U_1(x, 0^+)\), must satisfy the differential constraint \(2u'(x) = -f(x)\).
Therefore, since \(f \in L^2_T(\R)\) has zero mean, we uniquely reconstruct \(u\), ensuring uniqueness thanks to the zero-average condition \(\fint_0^T u(x) \ud x = 0\), proving \cref{eq: primitive_u}.

Since \(u\) satisfies the assumptions of \Cref{thm: elastic_dtn}, there exists a unique elastic extension \(\tilde{U}^+ \in H^1_T(\Hh^+_T, \R^2)\) that minimizes \(E_{\mathrm{el}}\) in the upper half-plane, yielding:
\begin{equation*}
	\frac{1}{2} \int_{\Hh^+_T} \sigma(\tilde{U}^+) : e(\tilde{U}^+) \ud x \ud y = \frac{G}{1-\nu} \frac{1}{4\pi} [u]^2_{\frac{1}{2},T}\,.
\end{equation*}
Moreover, \(\tilde{U}^+\) can be extended to the whole strip \(\Hh_T\) by enforcing the prescribed mirror symmetries. Therefore, we define the global minimizer as:
\begin{equation*}
	\tilde{U}(x,y) \coloneqq 
	\begin{cases}
		\tilde{U}^+(x,y), & \text{if } y > 0, \\
		\begin{pmatrix} -\tilde{U}_1^+(x,-y) \\ \tilde{U}_2^+(x,-y) \end{pmatrix}, & \text{if } y < 0.
	\end{cases}
\end{equation*}
Let us denote the restrictions to the half-planes \(\Hh^+_T\) and \(\Hh^-_T\), respectively, as \(\tilde{U}^+\) and \(\tilde{U}^-\). By construction, \(\tilde{U} \in SBV^2_T(\R^2, \R^2)\). Therefore, it easily follows that:
\begin{equation*}
	E_{\mathrm{el}}(\tilde{U}) = 2 \left( \frac{G}{1-\nu} \frac{1}{4\pi} [u]^2_{\frac{1}{2},T} \right) = \frac{G}{2\pi(1-\nu)} [u]^2_{\frac{1}{2},T}\,.
\end{equation*}
Setting \(c(G,\nu) \coloneqq \frac{G}{2\pi(1-\nu)}\), we obtain \cref{eq: l'en minima e' la seminorma}. 

We are left to observe that \(\tilde{U} \in AD(f)\). Indeed, the jump condition \(J(x) \coloneqq \tilde{U}^+(x,0) - \tilde{U}^-(x,0)=2u(x)\mathbf{e_1}\) is trivially satisfied and it is equivalent to \( \Curl \beta = f \mathbf{e_1} \mathcal{H}^1\llcorner\mathbb{H}^0_T\) in the sense of distributions. This concludes the proof.
\end{proof}
\begin{remark}
In the van der Merwe setting, the boundary traces of the global minimizer $U$ on the interface $\Hh^0_T$ are given by $U^\pm(x,0) = (\pm u(x), v(x))$, so that the jump across the interface is purely tangential, $[U](x) = 2u(x)\mathbf{e_1}$. Since \(-[U]' = f\mathbf{e_1}\) by \Cref{lemma: trace}, we obtain \(u' = -\frac{1}{2} f \in L^2_T(\R)\), hence \(u \in H^1_T(\R)\). Moreover, by \cref{eq: sol_generale2,eq: C_n relazione A_n e B_n,eq: A_n relazione B_n}, the normal component $v$ satisfies $\hat{v}_n = -i \operatorname{sgn}(\omega_n) \frac{1-2\nu}{2(1-\nu)} \hat{u}_n$; thus, also \(v \in H^1_T(\R)\). Consequently, \(U^\pm(\cdot, 0) \in H^1_T(\R, \R^2)\), resulting in a gain of half a derivative of regularity.
\end{remark}

\section{Read--Shockley energy regime} \label{sec: section4} \label{sec: read_shockley}
\subsection{Optimal energy regime}
From this subsection onward, we restrict ourselves to the setting relevant to the model under consideration. We assume that \(T\in\mathbb{N}\), where \(T\) is understood as the macroscopic periodicity scale of the problem, and we consider a dislocation density of the form
\begin{equation*}
    f=1-\sum_{i=1}^{T}\bigl(\delta_{x_i}\ast\rho^\theta\bigr),
\end{equation*}
corresponding to the presence of exactly \(T\) dislocations within each period. The assumption that there are precisely \(T\) dislocations in the interval \([0,T]\) does not bring any loss of generality. Indeed, once the periodicity scale and, hence, the number of dislocations contained in each period are fixed, one can rescale the space so that the length of the periodicity cell coincides with the number of dislocations.

Using this setting, we can take advantage of the results of \cite{Goldman26}.
Indeed, in that work, the authors prove that, in the different regimes depending on \(s\), the minimizer of the \((s,2)\)-Gagliardo seminorm for \(s \in (0,\frac{1}{2}]\) under the constraint in \cref{eq: constraint derivata} is the mantissa function \(u(x)=x-\lfloor x \rfloor\). Here, we refer to the critical case \(s=\frac{1}{2}\) where, in particular, the mantissa function has been regularized with a mollifier \(\rho^\theta =\frac{1}{\theta} \rho\left(\tfrac{\cdot}{\theta}\right)\) rescaled by \(\theta\), namely \(u^\theta(x) = u\ast \rho^\theta(x)\).
In the following result, we establish that the equispaced regularized configuration has an energetic logarithmic scaling.
\begin{proposition} \label{prop: blow up logaritmico}
Let \(u(x)=x-\lfloor x \rfloor\) and, for every \(0<\theta< \frac{1}{2}\), let \(u^\theta\) denote its mollification. Then
\[
[u^\theta]^2_{\frac{1}{2},T} = - 2 T \log \theta + O(1).
\]
\end{proposition}
\begin{proof}
Without loss of generality we can suppose the jump to be in \(x=0\) and, thanks to additivity and periodicity, we can take \(T=1\).
Proceeding as in \cite[Proposition 3.19]{Santilli26}, we compute 
\begin{align*}
	[u^\theta]^2_{\frac{1}{2},1} =& 2 \int_0^{+\infty} \frac{1}{h^2} \int_0^1 |u^\theta(x+h)-u^\theta(x)|^2 \ud x \ud h \\
	=&2 \int_0^{+\infty} \frac{1}{h^2} \int_{-\frac{1}{2}}^{\frac{1}{2}} |u^\theta(x+h)-u^\theta(x)|^2 \ud x \ud h \\
	=&2\int_0^{2\theta} \frac{1}{h^2} \int_{-\frac{1}{2}}^{\frac{1}{2}} |u^\theta(x+h)-u^\theta(x)|^2 \ud x \ud h + 2 \int_{2\theta}^{\frac{1}{2}}\frac{1}{h^2} \int_{-\frac{1}{2}}^{\frac{1}{2}} |u^\theta(x+h)-u^\theta(x)|^2 \ud x \ud h \\
	&+2\int_{\frac{1}{2}}^{+\infty} \frac{1}{h^2} \int_{-\frac{1}{2}}^{\frac{1}{2}} |u^\theta(x+h)-u^\theta(x)|^2 \ud x \ud h \\
	\eqqcolon&2(\mathcal{I}_1+\mathcal{I}_2 + \mathcal{I}_3),
\end{align*}
where in the first equality we have used the periodicity of the function. Now, we observe that \(\mathcal{I}_3=O(1)\) since 
\[
\mathcal{I}_3 \leq \int_{\frac{1}{2}}^{+\infty} \frac{1}{h^2} \ud h \leq 4.
\]
As for the first integral \(\mathcal{I}_1\), we split \(\left[-\frac{1}{2}, \frac{1}{2} \right]\) into two sets, which are 
\begin{align*}
E(h) &\coloneqq \left\{ x \in \left[-\frac{1}{2}, \frac{1}{2}\right] \suchthat [x,x+h] \cap [-\theta,\theta] = \varnothing \right\} \\
&= \left(\left[-\frac{1}{2}, \frac{1}{2}\right] \cap \left(-\infty,-\theta-h\right) \right)\cup \left(\left[-\frac{1}{2},\frac{1}{2}\right] \cap \left(\theta, +\infty\right)\right), \\
E^c(h) &\coloneqq \left[-\frac{1}{2}, \frac{1}{2} \right] \setminus E(h).
\end{align*}
Then, we point out that \(|E^c(h)|\leq|[-\theta,\theta]|+h\leq 4\theta\) and on this set, using the Mean Value Theorem, we easily get \(|u^\theta(x+h)-u^\theta(x)|\leq \frac{C(\|\rho'\|_\infty)h}{\theta}\). Outside of this set, using the same theorem we can easily get that \(|u^\theta(x+h)-u^\theta(x)|=h\). So, plugging in those estimates we get:
\begin{align*}
	\mathcal{I}_1&= \int_0^{2\theta} \frac{1}{h^2}\int_{E(h)} |u^\theta(x+h)-u^\theta(x)|^2 \ud x \ud h +\int_0^{2\theta} \frac{1}{h^2} \int_{E^c(h)} |u^\theta(x+h)-u^\theta(x)|^2 \ud x \ud h \\
	& \leq \int_0^{2\theta} \frac{1}{h^2} \left[ |E(h)| h^2+ |E^c(h)| \left( \frac{Ch}{\theta}\right)^2 \right] \ud h=O(1).
\end{align*}
We are now left with the second integral, where, again, we partition the domain into three sets:
\begin{align*}
	E_1(h)&\coloneqq\left\{x\in \left[-\frac{1}{2},\frac{1}{2}\right] \suchthat [x,x+h] \cap [-\theta,\theta] = \varnothing \right\},\\
	E_2(h)&\coloneqq \left\{x\in \left[-\frac{1}{2},\frac{1}{2}\right] \suchthat [-\theta,\theta] \subset [x,x+h] \right\},\\
	E_3(h)&\coloneqq \left\{x\in \left[-\frac{1}{2},\frac{1}{2}\right] \suchthat [-\theta,\theta] \cap [x,x+h] \neq \varnothing \right\} \setminus E_2(h),
\end{align*} 
from which we easily get
\begin{align*}
	\mathcal{I}_2&= \int_{2\theta}^{\frac{1}{2}} \frac{1}{h^2} \left[ \int_{E_1(h)} h^2 \ud x  + \int_{E_2(h)} (1-h)^2 \ud x +  \int_{E_3(h)}  |u^\theta(x+h)-u^\theta(x)|^2 \ud x \right] \ud h \\
	& =   \int_{2\theta}^{\frac{1}{2}} \frac{1}{h^2} \left(h^2 (1-2\theta-h)+ (h-2\theta)(1-h)^2 + O(\theta) \right) \ud h \\
	&= \int_{2\theta}^{\frac{1}{2}} \left( \frac{1}{h} - 1 + \frac{O(\theta)}{h^2} \right) \ud h= -\log\theta + O(1).
\end{align*}

\end{proof}

\subsection{Physical rescaling}
Our aim is now to reconstruct the geometric derivation proposed by Read and Shockley \cite{Read50} and to recover their energy regime from \Cref{prop: blow up logaritmico}. The authors see the atomic planes of the crystal as \virgolette{lines of flow} with a density of $\frac{1}{b}$, where $b = |\mathbf{b}|$ is the atomic distance. They consider a straight grain boundary oriented with an angle $\phi$ with respect to the \(y\)-axis, and each crystal is rotated by a small angle proportional to $\theta>0$. More precisely, the upper grain is rotated by $+\frac{\theta}{2}$ and the lower grain by $-\frac{\theta}{2}$. Consequently, the incident angles of the atomic planes intercepting the interface are geometrically defined as:
\begin{equation*}
\alpha \coloneqq \phi + \frac{\theta}{2} \quad \text{and} \quad \beta \coloneqq \phi - \frac{\theta}{2}\,.
\end{equation*}
By a standard goniometric argument, they deduce that the linear flux of planes matching the boundary from the upper grain is $\frac{1}{b} \cos \beta$, while the flux matching from the lower grain has density $\frac{1}{b} \cos \alpha$. The mismatch generated by this misorientation produces an edge dislocation. The linear density $\rho$ of these defects is thus given by the difference of the fluxes:
\begin{equation*}
\rho = \frac{1}{b}(\cos \beta - \cos \alpha)=\frac{2}{b} \sin(\phi) \sin\left(\frac{\theta}{2}\right)\,.
\end{equation*}
In the regime of small-angle grain boundaries ($\theta \ll 1$), we adopt the approximation \(\sin(x) \simeq x\) and, considering a symmetric tilt boundary ($\phi = \frac{\pi}{2}$), we obtain
\begin{equation*}
	\rho = \frac{\theta}{b}\,.
\end{equation*}

\begin{figure}[htpb]
	\centering
	\begin{tikzpicture}[scale=1.2, >=latex, font=\small]
		\definecolor{ottanio}{RGB}{0, 128, 128}
		\definecolor{ardesia}{RGB}{112, 128, 144}
		
		\def\bsp{1.0}        
		\def\halfTheta{22}   
		
		\pgfmathsetmacro{\dStep}{\bsp / sin(\halfTheta)}
		\pgfmathsetmacro{\angTop}{\halfTheta}
		\pgfmathsetmacro{\angBot}{-\halfTheta}
		
		\coordinate (A) at (0, 0);                 
		\coordinate (B) at (\dStep, 0);            
		\pgfmathsetmacro{\Cx}{\dStep - \bsp*sin(\halfTheta)}
		\pgfmathsetmacro{\Cy}{\bsp*cos(\halfTheta)}
		\coordinate (C) at (\Cx, \Cy);             
		
		\begin{scope}
			\clip (-4.8, -2.2) rectangle (4.8, 2.2);
			
			\foreach \k in {-4,...,2} {
				\draw[thick, ardesia!40] ({\k*\dStep}, 0) -- ++({\angTop}:18);
			}
			
			\foreach \k in {-4,...,2} {
				\draw[thick, ardesia!40] ({\k*\dStep}, 0) -- ++({\angBot}:18);
			}
		\end{scope}
		
		\draw[very thick, black] (-5.1, 0) node[left] {$\Hh^0_T$} -- (4.8, 0);
		
		\fill[ottanio!25] (A) -- (B) -- (C) -- cycle;
		\draw[very thick, ottanio] (A) -- (C);     
		\draw[very thick, ottanio] (A) -- (B);     
		\draw[very thick, violet] (B) -- (C);      
		
		\draw[thin, darkgray] (C) -- ++({\halfTheta+180}:0.18) -- ++({\halfTheta-90}:0.18) -- ++({\halfTheta}:0.18);
		
		\draw[<->, thick, violet] ([shift=({\halfTheta}:4pt)]C) -- ([shift=({\halfTheta}:4pt)]B)
			node[midway, above right=-2pt, font=\footnotesize\bfseries] {$b$};
			
		\draw[<->, thick, ottanio] ([yshift=-6pt]A) -- ([yshift=-6pt]B)
			node[midway, below=2pt, font=\footnotesize\bfseries, fill=white, inner sep=1.5pt] {$d = \frac{b}{\sin(\theta/2)}$};
			
		\draw[->, thick, ottanio] (1.3, 0) arc (0:\halfTheta:1.3);
		\node[font=\footnotesize, text=ottanio] at ({1.75*cos(\halfTheta/2)}, {1.75*sin(\halfTheta/2)}) {$\frac{\theta}{2}$};
	\end{tikzpicture}
	\caption{Read--Shockley construction for a symmetric tilt boundary ($\phi = \frac{\pi}{2}$). The right triangle highlights the atomic planes reaching the interface with spacing $d = \frac{b}{\sin(\theta/2)}$ and interplanar spacing $b$.}
	\label{fig:read_shockley_flux}
\end{figure}
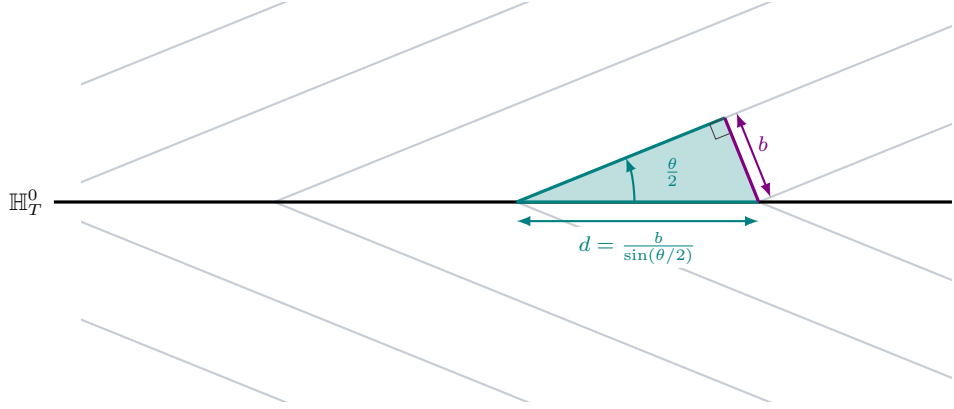

\begin{figure}[htpb]
\centering
\begin{tikzpicture}[
	scale=1.5,
	>=latex,
	font=\small
	]
	\coordinate (O) at (0,0);
	\coordinate (A) at (5,0);
	\coordinate (B) at (5,1.2);
	
	\draw[very thick, black] (O) -- (A) node[midway, below=4pt] {$\eps$ (Physical distance between dislocations)};
	
	\draw[very thick, violet, ->] (A) -- (B) node[midway, right=4pt] {$b$ (Burgers vector)};
	
	\draw[thick, dashed, darkgray] (O) -- (B) node[midway, above left] {Misoriented crystalline plane};
	
	\draw[thick, blue] (1.2,0) arc (0:13.5:1.2); 
	\node[blue] at (1.4, 0.15) {$\theta$};
	
	\draw[gray, thin] (4.8,0) -- (4.8,0.2) -- (5,0.2);
	
	\filldraw[black] (O) circle (1.5pt) node[left=6pt] {Previous defect};
	\filldraw[red] (B) circle (1.5pt) node[above=4pt] {Nucleation of new defect};
\end{tikzpicture}
\caption{Geometric derivation of the scaling identity $b = \eps \theta$.}
\label{fig:geometric_mismatch}
\end{figure}
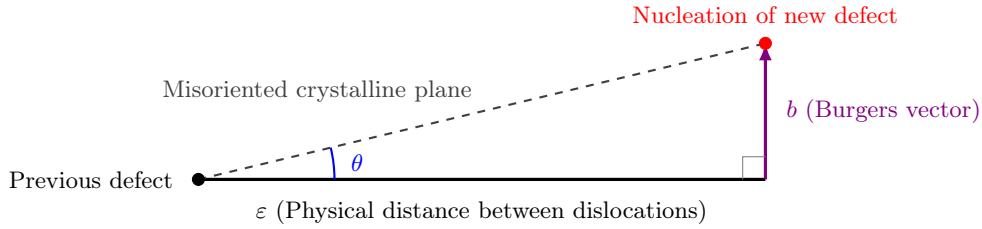

Given these relations, to evaluate the energy of the \emph{physical} optimal profile, it is necessary to establish the exact scaling between the physical dimensions of the crystal and the dimensionless variables of our normalized model. 
Let $\eps>0$ denote the distance between two adjacent dislocations along the interface. Since this microscopic spacing is the inverse of the linear density of the defects, we obtain	
\begin{equation} \label{eq: dislocation density}
	\eps = \frac{b}{\theta},
\end{equation}
yielding the relation \( b = \eps \theta \). To visualize this fundamental equivalence, one can refer to \cref{fig:geometric_mismatch}, where it is possible to recover the identity
\begin{equation*}
\tan(\theta) = \frac{b}{\eps}\,,
\end{equation*}
and, using the classical small-angle approximation \(\tan \theta \simeq \theta\), we recover \cref{eq: dislocation density}.

In the physical domain, the elastic singularity at the dislocation core is regularized over a microscopic distance comparable to the atomic spacing. Therefore, the physical mollification radius is proportional to the length of the Burgers vector, namely $r_{\text{phys}}\coloneqq cb $ for \(c>0\). From now on, without loss of generality, we will fix \(c=1\). In our mathematical framework, we study the regularized optimal profile on the normalized unit interval $[0,1]$. This implies a rescaling by \(\frac{1}{\eps}\) such that the distance between dislocations is normalized to $1$. Consequently, the mathematical mollification radius, which we use to regularize the periodic mantissa function, is evaluated as:
\begin{equation*}
r_{\text{math}} \coloneqq \frac{r_{\text{phys}}}{\eps} = \frac{\eps \theta}{\eps} = \theta\,.
\end{equation*}
This formally explains that, when computing the fractional seminorm on the normalized period (\Cref{prop: blow up logaritmico}), the correct mathematical mollification parameter to employ coincides with the misorientation angle $\theta$.
Let \(T_p \coloneqq T \eps\) denote the length of the physical interface containing \(T\) dislocations, and define 
\begin{equation*}
u^\theta_p(x) \coloneqq b u^\theta\left(\frac{x}{\eps}\right),
\end{equation*}
the rescaling of the regularized boundary profile. We define the physical elastic energy density of the interface as:
\begin{equation*}
	E_p \coloneqq \frac{C_{\mathrm{el}}}{T_p} [u^\theta_p]_{\frac{1}{2}, T_p}^2,
\end{equation*}
where \(C_{\mathrm{el}} \coloneqq \frac{c(G,\nu)}{4} = \frac{G}{8\pi(1-\nu)}\) takes into account the relation \(u' = -\frac{1}{2} f = -\frac{1}{2} (u^\theta)'\) and \Cref{thm: energia elastica e' la seminorma}.
\begin{corollary}
For every sufficiently small $\theta>0$, we have:
\begin{equation*}
	E_{p} = \frac{G b}{4\pi(1-\nu)} \, \theta |\log \theta| + O(\theta).
\end{equation*}
\end{corollary}

\begin{proof}
The physical fractional seminorm on the segment $[0, T_p]$ is given by:
\begin{align*}
	[u^\theta_p]_{\frac{1}{2}, T_p}^2 &= \int_0^{T_p} \int_\R \frac{|u^\theta_p(x+h) - u^\theta_p(x)|^2}{h^2} \ud h \ud x \\
	&=b^2 \int_0^T \int_\R \frac{\left|u^\theta(y+k) - u^\theta(y)\right|^2}{(\eps k)^2} (\eps \ud k) (\eps \ud y)\\
	&= b^2 [u^\theta]_{\frac{1}{2}, T}^2,
\end{align*}
where in the second equality we applied the change of variables $y = \frac{x}{\eps}$ and $k =\frac{h}{\eps}$. Exploiting the additivity of the energy over $T$ periods and applying \Cref{prop: blow up logaritmico} with mollification radius $\theta$, we obtain:
\begin{equation*}
	[u^\theta_p]_{\frac{1}{2}, T_p}^2 = b^2 T [u^\theta]_{\frac{1}{2}, 1}^2 = b^2 T \Big( -2 \log \theta + O(1) \Big)\,.
\end{equation*}
Finally, substituting this into the physical energy density and using \cref{eq: dislocation density}, we get:
\begin{align*}
	E_{p} &= \frac{C_{\mathrm{el}} [u^\theta_p]_{\frac{1}{2}, T_p}^2}{T_p} = C_{\mathrm{el}} \frac{b^2 T (-2 \log \theta + O(1))}{T \eps} \\
	&= \frac{G}{8\pi(1-\nu)} \frac{-2 b^2 \log \theta + O(1)}{b\theta^{-1}} = \frac{G b}{4\pi(1-\nu)} \, \theta |\log \theta| + O(\theta),
\end{align*}
concluding the proof.
\end{proof}
\begin{remark}
The same argument applies to the semi-coherent interface setting. Denoting by \(\delta>0\) the lattice-mismatch parameter, the corresponding regularized periodic profile satisfies
\[
    [u^\delta]_{\frac12,1}^2
    =
    2|\log\delta| + O(1)
    \qquad\text{as }\delta\to0.
\]
Thus, as in the Read--Shockley setting, the critical fractional energy exhibits a logarithmic divergence. After the physical rescaling, the geometric relation between the lattice mismatch and the spacing of the misfit dislocations introduces an additional factor of order \(\delta\), yielding
\[
    E_{\mathrm{VDM}}(\delta)
    =
    C\,\delta|\log\delta| + O(\delta),
\]
for a suitable elastic constant \(C>0\), and hence recovering the classical van der Merwe scaling.

A detailed derivation of this regime has already been carried out in \cite{DeLuca24}. In view of our dimensional-reduction result, their fractional formulation can be interpreted as the effective one-dimensional energy arising from the underlying two-dimensional elastic model. We therefore do not repeat the rescaling computation here.
\end{remark}

\section*{Acknowledgements}
The authors wish to express their deepest gratitude to Marcello Ponsiglione, Lucia De Luca and Emanuele Spadaro for their constant guidance and the insightful discussions throughout the preparation of this manuscript.
The authors are members of the Gruppo Nazionale per l'Analisi Matematica, la Probabilit\`a e le loro Applicazioni (GNAMPA) of the Istituto Nazionale di Alta Matematica (INdAM).

\sectionnotoc{AI Statement}
The authors acknowledge the use of AI tools in the preparation of this manuscript. These tools were used to revise the text and improve its exposition. The fundamental mathematical ideas and results are the work of the authors. The authors take full responsibility for the content of the manuscript.

\bibliographystyle{plainurl}
\bibliography{Bibliografia}

\par\vspace{2\baselineskip}
\small


(V. Corica) \textsc{Dip. di Matematica, Univ. Roma-I \virgolette{La Sapienza}, Piazzale Aldo Moro 5, 00185, Roma, Italy} \\
\textit{E-mail address, V. Corica:} \texttt{valerio.corica@uniroma1.it}

\par\medskip 

(F. Santilli) \textsc{Dip. di Matematica, Univ. Roma-I \virgolette{La Sapienza}, Piazzale Aldo Moro 5, 00185, Roma, Italy} \\
\textit{E-mail address, F. Santilli:} \texttt{francesco.santilli@uniroma1.it}


\end{document}